\documentclass{article}

\usepackage{xfrac}
\usepackage{amsfonts}
\usepackage{graphicx}
\usepackage{epstopdf}
\usepackage{algorithm}
\usepackage{algpseudocode}
\usepackage{amsfonts}
\usepackage{amsthm}
\usepackage{amsmath}
\usepackage{amssymb}
\usepackage{lipsum}
\usepackage{epstopdf}
\usepackage[utf8]{inputenc}
\usepackage{multirow}
\usepackage{authblk}
\usepackage{siunitx}
\usepackage{enumitem}
\usepackage{upgreek}
\usepackage{natbib}
\usepackage{subfigure}
\usepackage{url}

\usepackage{breakurl}
\usepackage[breaklinks]{hyperref}

\newtheorem{theorem}{Theorem}
\newtheorem{lemma}{Lemma}
\newtheorem{remark}{Remark}
\newtheorem{defn}{Definition}
\newtheorem{eg}{Example}

\title{Gradient ascent method for fully nonlinear parabolic differential equations with convex nonlinearity}

\author[1]{Hung Duong\thanks{\href{mailto:hduong1@utk.edu}{hduong1@utk.edu},. This work is the outcome of Hung Duong's PhD dissertation at the Florida State University.}}
\author[2]{Arash Fahim\thanks{\href{mailto:afahim@fsu.edu}{arash@math.fsu.edu}, \url{https://arashfahim.github.io}}}
\affil[1]{Department of Mathematics, University of Tennessee at Knoxville}
\affil[2]{Department of Mathematics, Florida State University}

\usepackage{amsopn}

\begin{document}

\maketitle

\begin{abstract}
We introduce a generic numerical schemes for fully nonlinear parabolic PDEs on the full domain, where the nonlinearity is convex on the Hessian of the solution. The main idea behind this paper is reduction of a fully nonlinear problem to a class of simpler semilinear ones parameterized by the diffusion term. The contribution of this paper is to provide a directional maximum principle with respect to the diffusion coefficient for semilinear problems, which specifies how to  modify the diffusion coefficient to approach to the solution of the fully nonlinear problem. While the objects of the study, diffusion coefficient, is infinite dimensional, the maximum direction of increase can be found explicitly.
This also provides a numerical gradient ascend method for the fully nonlinear problem. 
To establish a proof-of-concept, we test our method in a numerical experiment on the fully nonlinear Hamilton-Jacobi-Bellman equation for portfolio optimization under stochastic volatility model.
\end{abstract}

\maketitle


\section{Introduction}
\label{sec:intro}
Machine learning, specifically deep learning, has recently attracted considerable attention in the context of solving PDE numerically from a wide range of applications in physics, computer vision, and the like. See \citet{M24}. In this paper, we focus on a type of PDE, which appears in the area of continuous-time stochastic control, i.e., Hamilton-Jacobi-Bellman, HJB henceforth, equation. Such PDEs are often fully nonlinear. The more nonlinear a PDE gets, the harder it is to solve them. Further more, convergence and stability of numerical solutions become more challenging. Traditional methods of solving PDEs, such as finite-difference, and more contemporary ones, such as, Monte Carlo methods, are often slower on the nonlinearity. In this study, we specify three types of nonlinearity, which are  known in the stochastic control literature. A fully nonlinear \emph{parabolic} PDE is nonlinear on the Hessian of the solution with respect to $x$, i.e.,
\begin{equation}\label{eqn:fully_nonlinear}
 \begin{cases}
 	 	-\partial_t{v}-h(t,x,v(t,x),\nabla v(t,x),D^2v(t,x))=0& (t,x)\in Q\\
 	 	v(T,x)=g(x)&x\in \mathbb{R}^d
 \end{cases}
 \end{equation} 
where $Q:=[0,T)\times\mathbb{R}^d$, $\nabla v$ denotes the gradient of $v$ with respect to $x$, $D^2v$ denotes the Hessian of $u$ with respect to $x$,
$h:[0,T]\times Q\times \mathbb{R}\times \mathbb{R}^d\times \mathbb{M}_d\to\mathbb{R}\cup\{\infty\}$, where $\mathbb{M}_d$ is the set of all $d$ by $d$ matrices. 
Denote by $\mathcal{D}_h:=\{(t,x,p,q,\gamma): h(t,x,p,q,\gamma)<\infty\}$. Throughout this paper, we assume that the PDE is \emph{parabalic}, i.e., $h(t,x,p,q,\gamma)$ is nondecreasing in $\gamma$. Further more, we assume that $h(t,x,p,q,\gamma)$ is convex in $\gamma$.
An example of a fully nonlinear PDE is the Merton portfolio optimization problem, 
\begin{equation}\label{nonlinear_merton}
h(t,x,p,q,\gamma)=\begin{cases}
        -\dfrac{\mu^2p^2}{2\gamma}&\gamma<0\\
        \infty& \gamma\ge0
    \end{cases}
\end{equation}
A semilinear PDE is when the equation is nonlinear on the solution and/or its gradient with respect to $x$, i.e.,
\begin{equation}\label{eqn:semilinear}
 \begin{cases}
 	 	-\partial_t{v}-\frac12(\sigma^\intercal\sigma):D^2v(t,x)-h(t,x,v(t,x),\nabla v(t,x))=0& (t,x)\in Q\\
 	 	v(T,x)=g(x)&x\in \mathbb{R}^d
 \end{cases} 
\end{equation}
 In the above, for two square matrices $A$ and $B$, $A: B:=\textrm{Tr}[AB^\intercal]$ is the Frobenius  inner product and $B^\intercal$ is the transpose of matrix $B$. 
A linear PDE is when the equation is linear on all components but the independent variables,, i.e.,
\begin{equation*}
 \begin{cases}
 	 	-\partial_t{v}-\frac12(\sigma^\intercal\sigma):D^2v(t,x)-\mu\cdot\nabla v(t,x) + kv(t,x) - h(t,x)=0& (t,x)\in Q\\
 	 	v(T,x)=g(x)&x\in \mathbb{R}^d
 \end{cases} 
\end{equation*}

We first explain how to reduce a fully nonlinear PDE to a semilinear one.
Since $h(t,x,p,q,\gamma)$ is convex in $\gamma$, the Legendre-Fenchel transform of $h$ exists and is denoted by $f:(t,x,p,q,\sigma)\in[0,T]\times Q\times \mathbb{R}\times \mathbb{R}^d\times \mathbb{M}_d\to\mathbb{R}\cup\{\infty\}$:
 \begin{equation}
     \label{eqn:conjugate}
f(t,x,p,q,\sigma):=\sup_{\gamma\in \mathcal{D}_h(t,x,p,q)}\biggl\{\frac12(\sigma\sigma^\intercal): \gamma - h(t,x,p,q,\gamma)\biggl\}
 \end{equation}
Define $
    \mathcal{D}_f(t,x,p,q):= \{\sigma:f(t,x,p,q,\sigma)>-\infty\}$.  
    \begin{eg}\label{eg:simple_merton}
    For $h$ given in \eqref{nonlinear_merton}, $f(t,x,p,q,\sigma)=-|{\mu q }{\sigma}|$ with $\mathcal{D}_f(t,x,p,q) = \mathbb{R}$.
    \end{eg}
     We further 
 assume that $\mathcal{D}_f(t,x,p,q)$ is independent of $(t,x,p,q)$, and, therefore, we  drop $(t,x,p,q)$ and write $\mathcal{D}_f$.
By Fenchel–Moreau theorem,  
\begin{equation}
    \label{eqn:double_conjugate}  h(t,x,p,q,\gamma)=\sup_{\sigma\in \mathcal{D}_f}\biggl\{\frac12(\sigma^\intercal\sigma): \gamma - f(t,x,p,q,\sigma)\biggl\}
 \end{equation}
 \begin{remark}
The true notion of Legendre-Fenchel transform of $h$ is given by 
 \begin{equation}
f(t,x,p,q,a):=\sup_{\gamma\in \mathcal{D}_h(t,x,p,q)}\biggl\{\frac12a: \gamma - h(t,x,p,q,\gamma)\biggl\}
 \end{equation}
 where $a$ is an arbitrary square matrix, which makes $f$ convex in $a$. However, $\mathcal{D}_f(t,x,p,q)$ only includes  nonnegative-definite symmetric matrices. Hence, $a=\sigma^\intercal\sigma$ for some $\sigma\in \mathbb{M}_d$.
 \end{remark}
 \begin{remark}
     If $h$ is concave in $\gamma$, one can repeat the result of the paper with replacing $\sup$ with $\inf$ in \eqref{eqn:conjugate}. As a result, Section~\ref{sec:directional_maximum_principle} shall be titled ``Gradient descent''.
 \end{remark}
 
For a given $\sigma:[0,T]\times Q\to\mathbb{M}_d$ with $\sigma\in\mathcal{D}_f$, consider the following semi-linear problem
\begin{equation}\label{eqn:semi_nonlinear}
 \begin{cases}
 	 	-\partial_t{v}(t,x)-\frac12(\sigma^\intercal\sigma): D^2 v(t,x)+f(t,x,v(t,x),\nabla v(t,x),\sigma(t,x))=0\\
 	 	v(T,x)=g(x)
 \end{cases}
 \end{equation} 
It is well known that the solution to \eqref{eqn:semi_nonlinear} is related to the backward stochastic differential equations, BSDE henceforth, below:
\begin{equation}\label{eqn:bsde_L-F}
\begin{cases}
    dX^\sigma_t=\sigma(t,X^\sigma_t)dt\\
    Y^\sigma_t=g(X^\sigma_T)-\int_t^T\big(f(s,X^\sigma_s,Y^\sigma_s,Z^\sigma_s,\sigma_s)ds+Z^\sigma_s\sigma_sdB_s\big)
\end{cases}
\end{equation}
by $Y^\sigma_t=v(t,X^\sigma_t;\sigma)$ and $Z^\sigma_t=\nabla v(t,X^\sigma_t;\sigma)$, where $v(\cdot;\sigma)$ is the solution for \eqref{eqn:semi_nonlinear} and $B$ is a Brownian motion. In the first glance, the suggested relation requires differentiability of $v(\cdot,\sigma)$. However, the existence of the solution to the BSDE \eqref{eqn:bsde}, $(Y,Z))$, can be established independent of the semilinear PDE, which readily guarantees the existence of \emph{viscosity solution} for \eqref{eqn:semi_nonlinear}, without appealing to differentiability of $v(\cdot;\sigma)$. Moreover, the BSDE provides an efficient Monte Carlo numerical scheme for the semilinear problem \eqref{eqn:semi_nonlinear}.
For more information on the early work on BSDEs, see the book of \citet{MY99} and the references therein. The theory of viscosity solutions for nonlinear PDEs can be found in \citet{CIL92}. For Monte Carlo schemes based on BSDE see \citet{BT04} and \citet{Z01}.

The first attempt to generalize BSDE to apply to the fully nonlinear problems was \citet{CSTV07}, which introduces a second order BSDE,  2BSDE henceforth, of the form 
\begin{equation}\label{eqn:2BSDE_CSTV}
\begin{cases}
    dX_t=\sigma(t,X_t)dt\\
Y_t=g(X_T)+\int_t^T\big(h_0(s,X_s,Y_s,Z_s,\Gamma_s)ds-Z_s\sigma_sdB_s\big)\\
    Z_t=\nabla g(X_T)-\int_t^T \Gamma_s\sigma_sdB_s
    \end{cases}
\end{equation}
to represent the fully nonlinear PDEs of the form
\begin{equation}\label{eqn:linear+fullynonlinear}
 \begin{cases}
 	 	-\partial_t{v}(t,x)-\frac12(\sigma^\intercal\sigma): D^2v(t,x)-h_0(t,x,v(t,x),\nabla v(t,x)&\!\!\!\!\!\!,D^2v(t,x))=0 \\&(t,x)\in Q\\
 	 	v(T,x)=g(x)&x\in \mathbb{R}^d
 \end{cases}
 \end{equation}
 Note that while $h(t,x,p,q,\gamma)=\frac12(\sigma^\intercal\sigma): \gamma+h_0(t,x,p,q,\gamma)$ is fully nonlinear, it requires to have a linear component, $\frac12(\sigma^\intercal\sigma): \gamma$. Similar to the BSDE, the solution of the 2BSDE, $(Y,Z,\Gamma)$, are related to the solution of \eqref{eqn:linear+fullynonlinear} through $Y_t=v(t,X_t)$, $Z_t=\nabla v(t,X_t)$, and $\Gamma_t=D^2v(t,X_t)$. However, \citet{CSTV07} only showed the uniqueness of the solution $(Y,Z,\Gamma)$ and didn't provide any existence except when the PDE has a strong solution. Nevertheless, some studies, \citet{FTW11}, \citet{T13}, \citet{AA13}, and \citet{GZZ15}, used the 2BSDE in \citet{CSTV07} to propose Monte Carlo schemes for \eqref{eqn:fully_nonlinear} by adding a small viscosity term to the  fully nonlinear problem, when the linear term is missing. 
 
In \cite{STZ12}, they leveraged convexity of $h$ in $\gamma$ and Legendre-Fenchel transform \eqref{eqn:conjugate} to propose a new theory for 2BSDE with proper existence and uniqueness results, which we review in more details in Section~\ref{sec:max_principle}. This closes the gap in the earlier paper \citet{CSTV07}.
More specifically, \citet{CSTV07} shows that if \eqref{eqn:linear+fullynonlinear} has a smooth solution, then $Y_t=v(t,X_t)$, $Z_t=\nabla v(t,X_t)$, and $\Gamma_t=D^2v(t,X_t)$ satisfy \eqref{eqn:2BSDE_CSTV}. Apart form the assumption of regularity of the solution of the fully nonlinear PDE, the result of \citet{CSTV07} requires the existence of a linear term in the fully nonlinear PDE  \eqref{eqn:2BSDE_CSTV} alongside the nonlinearity $h_0$. The reliance on the linear term is relaxed in the work of \citet{STZ12}. 

To present a heuristic interpretation of the result of  \cite{STZ12} in the Markovian case, we first denote by $v(t,x;\sigma)$ the solution for \eqref{eqn:semi_nonlinear}, with diffusion coefficient $\sigma$. Then,  \eqref{eqn:conjugate} implies that
\begin{equation}   
h(t,x,p,q,\gamma) \ge\frac12(\sigma^\intercal\sigma): \gamma - f(t,x,p,q,\sigma)
\end{equation}
and, therefore,
$v(t,x;\sigma)$ constitute a subsolution for \eqref{eqn:fully_nonlinear}. Given the comparison principle holds for fully nonlinear PDE within a suitable class of functions, we have $v(t,x)\ge v(t,x;\sigma)$, and therefore, $v(t,x)\ge \sup_{\sigma\in\mathcal{D}_f}v(t,x;\sigma)$.
The result of \cite{STZ12} guarantees that the inequality is an equality:
\begin{equation}\label{sup:semi-linear}
v(t,x)=\sup_{\sigma\in\mathcal{D}_f}v(t,x;\sigma) 
\end{equation} 
It is important to note that the right-hand side of \eqref{sup:semi-linear} is not guaranteed to exists or be well-defined, unless we have sufficient regularity properties for $v(t,x;\sigma)$ for any $\sigma\in\mathcal{D}_f$.  

The result of \citet{STZ12} is analogous to Perron's method for viscosity solutions, which obtains a solution by taking the supremum of subsolutions (sub harmonic functions). Here, \eqref{sup:semi-linear} takes a supremum of a specific class of subsolutions of \eqref{eqn:fully_nonlinear} to obtain a solution to \eqref{eqn:fully_nonlinear}.

The main contribution of this paper is to present a systematic gradient ascent approach to adjust $\sigma$ such that   $v(t,x;\sigma)$ increases toward $v(t,x)$, the solution of \eqref{eqn:fully_nonlinear}. The classical maximum principle for semilinear PDEs, equivalently BSDEs, focuses on increasing a parameterized source term, $f(t,x,p,q,\theta)$ with parameter $\theta$ independent of $\sigma$, while keeping $\sigma$ constant. In other words, we have  $\sup_\theta v(t,x;\theta)=v(t,x)$, where $v(t,x;\theta)$ and $v(t,x)$ are solutions to \eqref{eqn:semilinear} with source term $h(t,x,p,q) = -f(t,x,p,q,\theta)$ and $h(t,x,p,q) = -\inf_\theta f(t,x,p,q,\theta)$, respectively.
In the case of \eqref{eqn:semi_nonlinear},  the parameter of the source term is the diffusion coefficient itself, $\sigma$. The new maximum principle is based on a notion of directional derivative of $v(t,x;\sigma)$ in $\sigma(t,x)$,
\begin{equation}
    {\nabla}_{\varsigma} v(t,x;\sigma):=\lim_{\epsilon\to 0}\frac{v(t,x;\sigma+\epsilon \varsigma)-v(t,x;\sigma)}{\epsilon},~~~ 
\end{equation}
where $\varsigma:[0,T]\times Q\to\mathbb{M}_d$ and $\sigma, \sigma+\epsilon \varsigma\in\textrm{\normalfont relint}(\mathcal{D}_f)$, for all $\epsilon\in[0,\epsilon_0]$ for some $\epsilon_0>$. This directional derivative is well-defined if $v(t,x;\sigma)$ satisfies certain continuous dependence result on $\sigma$. For the purpose of this paper, we assume that such continuous dependence holds and focus of evaluation of the directional derivative ${\nabla}_{\varsigma} v(t,x;\sigma)$. For the Hamilton-Jacobi-Bellman equations from stochastic control problems, continuous dependence result is established, for instance, in \citet{JK2002}.

While we show that For any $\varsigma$, Theorem~\ref{thm:max_principle} shows that 
${\nabla}_{\varsigma} v(t,x;\sigma)$ 
satisfies a linear partial differential equation with coefficients depending on $\varsigma(t,x)$, $v(t,x;\sigma)$, $\nabla v(t,x;\sigma)$,  and $D^2v(t,x;\sigma)$. Based on this result, Theorem~\ref{thm:positive_gradient} provides the direction of maximum increase of $v(t,x;\sigma)$, $\varsigma^*$ in terms of $\sigma$, $v(t,x;\sigma)$, $\nabla v(t,x;\sigma)$, and $D^2v(t,x;\sigma)$. This allows us to modify $\sigma$ to $\sigma+\alpha_n{\nabla}_{\varsigma*} v(t,x;\sigma)\varsigma^*$ for some $\alpha_n>0$ in the gradient ascent. Finally, Theorem~\ref{thm:derive_zero}, \ref{thm:derive_zero} asserts the condition under which the maximizer of $v(t,x;\sigma)$ on $\sigma$ is attained.  

To make the gradient ascent method practical, we shall use an scheme which yields an approximate solution to linear and semilinear PDEs, which yields a \emph{field estimate} of the solution and its derivatives. A field estimate shall be a function rather than a set of value in discrete points. Therefore, some early numerical methods for semilinear PDEs  such as \citet{BT04} and \citet{GLW05},require an extra interpolation step to obtain the solution. Moreover, the interpolation step must provide an accurate interpolation for the first and second derivative simultaneously, which is cumbersome and sometimes impossible. Therefore, the prefered method of solving semilinear PDEs is the deep numerical schemes such as the one in \citet{JHW18} and in the preprint version \citet{JHW17}. In our study, we specifically use the deep scheme in \citet{JHW17} and \citet{JHW18} to approximate semilinear PDEs, because such approximations readily give us a global estimate of the solution and its first derivatives. Even though there is no guarantee that the deep scheme in \citet{JHW17} and \citet{JHW18} provides an accurate estimation of the second derivatives, we shall see in the numerical experiment that it does not pose a major challenge and using more computational resources can lead to an accurate second derivative. While the main result of the paper is theoretical, for the purpose of completeness, in Appendix~\ref{sec:semi_linear_scheme}, we provide a short description of the deep learning scheme of \citet{JHW18}. 

For fully nonlinear problems, deep numerical schemes are introduced by \citet{BWJ19}, \citet{GPW21}, and \citet{GPW22}. Such schemes are approaching the problem from the angle of \citet{CSTV07} and work best for PDEs of the form \eqref{eqn:linear+fullynonlinear}. In the current paper, we do not intend to establish a comparison between the two approaches. However, our implementations of both algorithms show that \citet{BWJ19}, \citet{GPW21}, and \citet{GPW22} provide much faster schemes in the problems where the fully nonlinear problem has a linear term. Therefore, we only recommend our approach to problems which are not of the form \eqref{eqn:linear+fullynonlinear}.

The paper is organized as follows. In Section~\ref{sec:max_principle}, we briefly present the relevant results from 2BSDEs and the maximum principle, which we use to establish the proposed numerical scheme. All the theoretical  results of the paper are presented in this section. In Section~\ref{sec:numerics}, we test our scheme on the fully nonlinear PDE coming from a stochastic volatility model in quantitative finance. We postpone all the proofs to the appendix.

\section{Preliminaries}
\label{sec:prelim}
The connection between fully nonlinear convex parabolic PDEs  and 2BSDEs was first studied in \citet{CSTV07} for the case when the PDE has a classical solution. The rigorous study of 2BSDEs and their characterization in terms of the maximum of a class of BSDEs is later established in \citet[Theorem~5.3]{STZ12}. To set the base for our main results, we provide a summary of the latter paper. 

Let $X$ be the canonical process on a suitable extension of Wiener space, $\Omega:=\textrm{\normalfont C}([0,T];\mathbb{R}^d)$ with an augmented right-continuous filtration and $\langle X\rangle$ represent the quadratic variation of process $X$ as defined pathwise by \citet{Karandikar}. 
We denote the Wiener measure on $\Omega$ by $\mathbb{P}_0$.For a progressively measurable process $\sigma$ with values in $\mathbb{M}_d$, $\mathbb{P}^\sigma$ is given by $\mathbb{P}^\sigma=\Pi^{-1}\circ \mathbb{P}_0s$, where $\Pi:X\mapsto\int_0^\cdot \sigma^{-1}_sdX_s$. 

The 2BSDE associated to \eqref{eqn:fully_nonlinear} is given by
\begin{equation}\label{eqn:2bsde}
    Y_t=g(X_T)-\int_t^Tf(s,X_s,Y_s,Z_s,\sigma_s)ds-\int_t^T Z_s dX_s+K^\sigma_T-K^\sigma_t,~~\mathbb{P}^\sigma\textrm{-a.s.} 
\end{equation}
where $\mathbb{P}^\sigma$ is a probability measure under which $d\langle X\rangle_t=\sigma_t^\intercal\sigma_t  dt$ and $K$ is a nondecreasing process $\mathbb{P}^\sigma$-a.s. 
A solution to \eqref{eqn:2bsde} is defined to be a pair of processes $(Y,Z)$ such that for any progressively measurable process $\sigma$ with values in $\mathcal{D}_f$, 
\begin{equation}
    K_t^\sigma:=Y_0-Y_t+\int_0^tf(s,X_s,Y_s,Z_s,\sigma_s)ds+\int_0^t Z_s dX_s,~~\mathbb{P}^\sigma\textrm{-a.s.}
\end{equation}
is nondecreaing $\mathbb{P}^\sigma$-a.s. and satisfies the minimality condition
\begin{equation}\label{cond:minimality_K}
    K_t^\sigma= \textrm{\normalfont essinf}_{\sigma\in S_t^\sigma}\mathbb{E}^{\sigma}[K^{\sigma}_1], ~~\mathbb{P}^\sigma\textrm{-a.s. for all }
\end{equation}
where $S_t^\sigma$ is the set of all progressively measurable processes $\sigma$ with values in $\mathbb{M}_d$ such that $\sigma(s,\omega)=\sigma(s,\omega)$ on $(s,\omega)\in[0,t]\times\Omega$. Specifically, the minimality condition \eqref{cond:minimality_K} ensures that $K$ can be universally defined and is equal to $K^\sigma$ on the support of $\mathbb{P}^\sigma$-a.s. 
\begin{remark}\label{rem:mutual_singularity}
Let $\sigma$ be a progressively measurable process  with values in $\mathbb{M}_d$ such that $\int_0^T\sigma_s^\intercal\sigma_sds<\infty$.
By $\mathcal{P}$, we denote the set of probability measures $\mathbb{P}$ on $\Omega$ such that $\mathbb{P}=\Pi^{-1}\circ \mathbb{P}_0$, where $\Pi:X\mapsto\int_0^\cdot \sigma^{-1}_sdX_s$. $\mathbb{P}_0$-a.s.
    The 2BSDE \eqref{eqn:2bsde} is held $\mathbb{P}$-a.s. for all $\mathbb{P}\in\mathcal{P}$. This is referred to as $\mathbb{P}$-quasi-surely or $\mathbb{P}$-q.s. equality. Since for $\sigma\neq\sigma^\prime$, $\mathbb{P}^\sigma$ and $\mathbb{P}^{\sigma^\prime}$ are mutually singular, i.e., the support of $\mathbb{P}^\sigma$ and $\mathbb{P}^{\sigma^\prime}$ are mutually exclusive, \eqref{eqn:2bsde} is valid on mutually exclusive slices of the Wiener space $\Omega$. 
\end{remark}
For $p\ge0$, we say a function $\psi$ satisfies $p^\textrm{th}$ growth condition if $|\psi(t,x)|\le C(1+|x|^p)$. By $\textrm{\normalfont C}_p^{1,2}$, we denote the space of all functions $\phi(t,x):[0,T]\times Q\to\mathbb{R}$ such that  $\partial_t \phi$, $\nabla \phi$, and  $D^2 \phi$ exist and are continuous in $(t,x)$ and satisfy the $p^\textrm{th}$ growth condition.
If a solution $v$ to \eqref{eqn:fully_nonlinear} belongs to $\textrm{\normalfont C}_p^{1,2}$, for any $\sigma\in\mathcal{D}_f$ such that $dX_t=\sigma_tdB_t$ has a weak solution with $\mathbb{E}[\sup_{t\in[0,T]}|X_t|^p]<\infty$, we have $Y_t=v(t,X_t)$, $Z_t=\nabla v(t,X_t)$, and  $\gamma_t=D^2v(t,X_t)$ satisfy \eqref{eqn:bsde} with $K^\sigma$ given by
\begin{equation}\label{cond:minimality_K_Markov}
    K_t^\sigma:=\int_0^t\Big(
    h(s,X_s,Y_s,Z_s,\gamma_s)-\frac12(\sigma^\intercal_s\sigma_s): \gamma_s+f(s,X_s,Y_s,Z_s,\sigma_s)
    \Big)ds
\end{equation}
is nondecreasing, because \eqref{eqn:conjugate} implies that for any $\sigma$, 
\begin{equation}
    h(s,X_s,Y_s,Z_s,\gamma_s)-\frac12(\sigma^\intercal_s\sigma_s): \gamma_s+f(s,X_s,Y_s,Z_s,\sigma_s)\ge 0
\end{equation}
Note that in \eqref{cond:minimality_K_Markov}, $d\langle X\rangle_t=\sigma^\intercal_t\sigma_tdt$, which may not be necessarily Markovian. 
However, if there exists 
\begin{equation}\label{eqn:minimizing_K}
    \begin{split}
        \sigma^*(t,x)\in\mathop{\textrm{\normalfont argmin}}\limits_{\sigma\in\mathcal{D}_f}\Big\{&h\big(t,x,v(t,x),\nabla v(t,x),D^2v(t,x))\\
        &-\frac12(\sigma^{\intercal}\sigma): D^2 v(t,x)+f\big(t,x,v(t,x),\nabla v(t,x),\sigma(t,x)\big)\Big\}
x    \end{split}
\end{equation}
such that $dX_t=\sigma^*(t,X_t)dB_t$ has a weak solution with $\mathbb{E}[\sup_{t\in[0,T]}|X_t|^p]<\infty$, then, $\mathbb{P}^{\sigma^*}$-a.s.,   we have $K^{\sigma^*}\equiv0$  and 
\begin{equation}
\begin{split}
        Y_t&=g(X_T)-\int_t^Tf(s,X_s,Y_s,Z_s,\sigma^*_s)ds-\int_t^T Z_s dX_s\\
    dX_t&=\sigma^*(t,X_t)dB_t
\end{split}
\end{equation}
\begin{eg}\label{eg:simple_merton_min_K}
    For $h$ given in \eqref{nonlinear_merton}, Example~\ref{eg:simple_merton} provides $f(t,x,p,q,\sigma)=-|{\mu q }{\sigma}|$ and, therefore, \eqref{eqn:minimizing_K} suggests $\sigma^*(t,x)=-\frac{|{\mu v_x }|}{v_{xx}}$, where $v(t,x)$ is the solution to the fully nonlinear equation.
\end{eg}

\subsection{Maximum principle for 2BSDEs}
\label{sec:max_principle_deep}
Another approach to solve \eqref{eqn:fully_nonlinear} is the maximum principle approach for 2BSDE. Specifically, \citet[Theorem~4.3]{STZ12} asserts that 
\begin{equation}\label{eqn:sup_bsde}
    Y_t=\textrm{\normalfont esssup}_{\sigma\in S_t} Y^\sigma_t
\end{equation}
where $Y^\sigma_t$ satisfies the BSDE
\begin{equation}
    \label{eqn:bsde}
    Y^\sigma_t=g(X_T)-\int_t^T f(s,X_s,Y^\sigma_s,Z^\sigma_s,\sigma_s)ds-\int_t^T Z^\sigma_s dX_s,~~~dX_s=\sigma_sdB_s
\end{equation}
In the Markovian cases, if we restrict the supremum in \eqref{eqn:sup_bsde} to Markovian $\sigma$, then $Y^\sigma_t=v(t,X_t;\sigma)$ where $v(\cdot;\sigma)$ is the solution to \eqref{eqn:semi_nonlinear}. Therefore, \eqref{eqn:sup_bsde} can be written as 
\begin{equation}
    Y_t=\textrm{\normalfont esssup}_{\sigma\in S_t} v(t,X_t;\sigma)
\end{equation}
Recall from previous section that $Y_t=v(t,X_t)$, where $v$ is the solution to \eqref{eqn:fully_nonlinear} and $dX_t=\sigma^*(t,X_t)dB_t$ with $\sigma^*$ given by \eqref{eqn:minimizing_K}.  
\begin{equation}\label{eqn:maximizer_sigma}
    v(t,X_t)=v(t,X_t;\sigma^*) =\textrm{\normalfont esssup}_{\sigma\in S_t} v(t,X_t;\sigma)
\end{equation}
In the above, $v(t,x;\sigma)$ is the solution to \eqref{eqn:semi_nonlinear}.
\begin{eg}
    Following up on Example~\ref{eg:simple_merton_min_K},  the solution the solution to \eqref{eqn:semi_nonlinear} with
    $f(t,x,p,q,\sigma)=-|{\mu q }{\sigma}|$  for $\sigma$ constant can be evaluate in closed from when $g(x)=1-e^{-\eta x}$: 
    \begin{equation}
        v(t,x;\sigma)=1-e^{-\eta x+\alpha(T-t)},~\textrm{ with }~\alpha=\frac12\sigma^2\eta^2-|\lambda|\sigma\eta
    \end{equation}
    Then, $\sup_{\sigma\in\mathbb{R}}v(t,x;\sigma)$ wields $\sigma^*(t,x)=\frac{|\lambda|}{\eta}$. In this case,  $v(t,x;\sigma)=1-e^{-\eta x+\alpha^*(T-t)}$ with $\alpha^* = -\frac{|\lambda|}{2}$ satisfies the fully nonlinear problem \eqref{eqn:fully_nonlinear}. However, we should emphasize that in general the supremum has to be over all diffusion coefficients and not only constabt diffusions.
\end{eg}

\section{Main Results}
\label{sec:max_principle}

\subsection{Directional maximum principle}
\label{sec:directional_maximum_principle}
The main challenge in implementing \eqref{eqn:maximizer_sigma} is to find how to modify $\sigma\in S_t$ to increase the value of $v(t,x;\sigma)$ for all $(t,x)$. 
The idea behind the direction maximum principle is to search  for a sequence $\{\sigma_n\}_{n\ge1} \subset S_t$ such that $\sigma_{n+1}-\sigma_n$ indicates the maximum direction of increase, if it exists, for $v(t,x;\sigma_n)$ for all $(t,x)$. Then, we hope to establish
\begin{equation}
    \lim_{n\to\infty} v(t,x;\sigma_n)=v(t,x;\sigma^*),
\end{equation}
To obtain the maximum direction of increase, we appeal to gradient ascent procedure. To define the ``derivative of $v(\cdot;\sigma)$ with respect to $\sigma$'' appropriate to our study, we appeal to the notion of directional derivative. For a $\varsigma:[0,T]\times Q\to\mathbb{M}_d$, we loosely define
\begin{equation}\label{denf:directional_derivative}
    {\nabla}_{\varsigma} v(\cdot;\sigma):=\lim_{\epsilon\to 0}\frac{v(\cdot;\sigma+\epsilon \varsigma)-v(\cdot;\sigma)}{\epsilon},~~~ 
\end{equation}
Theorem~\ref{thm:max_principle} below shows that if ${\nabla}_{\varsigma} v(\cdot;\sigma)$ is well defined, then it satisfies the linear parabolic PDE
\begin{equation}\label{eqn: directional derivative of u}
\begin{cases}
-\partial_{t}u-\frac{1}{2}\sigma^\intercal\sigma \!: \! D^{2}u  \!+\! ku - \mu\!\cdot\! \nabla u+\ell:\varsigma=0
\\ 
u(T,x)=0
\end{cases}
\end{equation}
where 
\begin{equation}\label{eqn:linear_coefficients}
    \begin{split}
        k(t,x)&:=\partial_p f\big(t,x,v(t,x;\sigma),\nabla v(t,x;\sigma),\sigma(t,x)\big)\\
        \mu(t,x)&:=-
\partial_{q}f\big(t,x,v,\nabla v(;\sigma),\sigma\big)\\
        \ell(t,x)&:=\partial_{\sigma}f\big(t,x,v(t,x;\sigma),\nabla v(t,x;\sigma),\sigma(t,x)\big)-\sigma^\intercal  D^2v(t,x;\sigma)
    \end{split}
\end{equation}
In the above, $\partial_{p}f$, $\partial_{q}f$, and $\partial_{\sigma}f$ are partial derivatives of $f$ with respect to $p$, gradient of $f$ with respect to $q$, and gradient of $f$ with respect to $\sigma$, respectively, and for two vectors $\mu$ and $\nabla u$, $\mu\cdot\nabla u$ denotes the inner product. In the following, we provide a necessary condition, for \eqref{eqn: directional derivative of u} to exits and for Theorems \ref{thm:max_principle} and \ref{thm:positive_gradient}, which provide us with a theoretical representation of the direction of maximum increase and optimality condition for $\sigma$.

We recall that $\textrm{\normalfont relint}(A)$ denotes the relative interior of a set $A$ of functions equipped with $L^\infty$-norm. 
In the following results, we assume that $\sigma$ satisfies conditions in regard to the regularity and continuous dependence properties of the solution to the semilinear equation \eqref{eqn:semi_nonlinear}.
\begin{defn}\label{def:regularity}
We say $\sigma:[0,T]\times Q\to\mathbb{M}_d$ with
     $\sigma\in\mathcal{D}_f$ satisfies \emph{regularity assumption} if 
    \begin{enumerate}
    \item $\sigma^{-1}$ is bounded and $\int_0^T\sigma(t,\omega(t))^\intercal\sigma(t,\omega(t))dt<\infty$, where $t\to\omega(t)$ is the canonical process in the Wiener process.
        \item the semilinear \eqref{eqn:semi_nonlinear} admits a solution, $v(\cdot;\sigma)$, in $\textrm{\normalfont C}_p^{1,2}$, for some $p\ge0$, and $\nabla v(\cdot;\sigma)$ and $D^2v(\cdot;\sigma)$ are Lipschitz on $(t,x)$.
        \item the linear equation \eqref{eqn: directional derivative of u} satisfies comparison principle for viscosity solutions in the class of functions with growth $p\ge0$, i.e., $|\phi(t,x)|\le C(1+|x|^p)$.
    \end{enumerate}
\end{defn}
\begin{remark}\label{rem:Girsanov}
Let $X$ be the canonical process in the Wiener space. Recall from Remark~\ref{rem:mutual_singularity} that
for a given $\sigma:[0,T]\times Q\to\mathbb{M}_d$ such that $\int_0^T\sigma(t,X_t)^\intercal\sigma(t,X_t)dt<\infty$. $\Pi:X\mapsto \int_0^\cdot\sigma^{-1}(t,X_t)dX_t$ defined a map from the Wiener space to itself and $\mathbb{P}=\Pi^{-1}\circ \mathbb{P}_0$ provides a weak solution for $dX_t=\sigma(t,X_t) dB_t$. 
     Then, by (1) in Definition~\ref{def:regularity},  Girsanov change of measure provides a solution for $dX_s=\mu(s,X_s)ds+\sigma(s,X_s)dB_s$ for any bounded $\mu:[0,T]\times Q\to\mathbb{R}^d$.
\end{remark}
\begin{defn}\label{def:cont_dep}
    We say $\sigma\in\mathcal{D}_f$ satisfies \emph{continuous dependence} in the direction of $\varsigma\in\mathbb{M}_d$, if  for some $\epsilon_0>0$, and  all $\epsilon\in(0,\epsilon_0)$, $\sigma+\epsilon\varsigma\in{\textrm{\normalfont relint}(\mathcal{D}_f)}$, \eqref{eqn:semi_nonlinear} admits a continuous viscosity solution $v(\cdot;\sigma+\epsilon\varsigma)$, and there exists a bounded function $\kappa:[0,T]\to\mathbb{R}_+$ and $p\ge0$ such that for some 
    \begin{equation}\label{eqn:continuous_dependency}
           |v(t,x;\sigma)-v(t,x;\sigma+\epsilon\varsigma)|\le \epsilon (1+|x|^p)\kappa(T-t)
\end{equation}
with $\lim_{t\to0}\kappa(t)=0$, and $\kappa$ and $p$ are independent of $\epsilon\in(0,\epsilon_0)$ and 
\end{defn}
We impose the following assumptions in this section:

\begin{enumerate}[label =   \textbf{A\arabic*.}]
\item Problem \eqref{eqn:fully_nonlinear} admits comparison principle in the class of functions with growth $p\ge0$.
    \item $f:[0,T]\times Q\times\mathbb{R}\times\mathbb{R}^d\times\mathbb{M}_d\to\mathbb{R}$ is continuously differentiable on $(t,x,p,q,\sigma)$ $\in \mathcal{D}_f$ with bounded Lipschitz continuous derivatives, $\mathcal{D}_f$ is independent of $(t,x,p,q)$, and  $\mathcal{D}_f={\textrm{\normalfont relint}(\mathcal{D}_f)}$
    \item 
 $\sigma\in\mathcal{D}_f$ satisfies regularity and continuous dependence in direction of $\varsigma$, i.e., Definitions~\ref{def:regularity}-\ref{def:cont_dep}.
\end{enumerate}
\begin{remark}\label{rem:A2}
    Assumption~\textbf{A2} ensures that for a fixed $\sigma$ such that $(t,x,v(t,x;\sigma), $ $\nabla v(t,x;\sigma),$ $\sigma(t,x))\in\mathcal{D}_f$ and $\sigma$ satisfies continuous dependence in the direction of $\varsigma$, then $(t,x,v(t,x;\sigma+\epsilon\varsigma),\nabla v(t,x;\sigma+\epsilon\varsigma),(\sigma+\epsilon\varsigma)(t,x))\in\mathcal{D}_f$ independent of $(t,x)$. 
\end{remark}
\begin{theorem}\label{thm:max_principle} 
Under assumptions \textbf{A2-A3}, ${\nabla}_{\varsigma} v(\cdot;\sigma)$ is well-defined and is a continuous viscosity solution to linear parabolic problem \eqref{eqn: directional derivative of u}.
\end{theorem}

The promise of gradient ascent is to find $\varsigma$ such that $\nabla_\varsigma v(\cdot;\sigma)>0$, desirably uniformly in $(t,x)$. This implies that for some sufficiently small  $\delta, \epsilon>0$,   
\begin{equation}
    \dfrac{v(\cdot;\sigma+\epsilon\varsigma)-v(\cdot;\sigma)}{\epsilon}\ge \nabla_\varsigma v(\cdot;\sigma)-\delta>0
\end{equation}
which implies that $v(\cdot;\sigma+\epsilon\varsigma)>v(\cdot;\sigma)$. The following Theorem provides a sufficient condition for $\nabla_\varsigma v(\cdot;\sigma)>0$. 
\begin{theorem}\label{thm:positive_gradient}
Let Assumptions \textbf{A2-A3} for $\varsigma=-\alpha\ell$, where $\ell$ is given in \eqref{eqn:linear_coefficients} and $\alpha:[0,T]\times\mathbb{R}^2\to\mathbb{R}$ is a positive function with $\alpha<\epsilon_0$, where $\epsilon_0$ is as in Definition~\ref{def:cont_dep}. 
Further  assume that $\ell\not\equiv0$. Then,   $\nabla_\varsigma v(\cdot;\sigma)>0$.
\end{theorem}

The above theorem, guarantees a choice of $\varsigma$ such that $\nabla_\varsigma v(\cdot;\sigma)>0$ and, therefore,  $v(\cdot;\sigma+\epsilon\varsigma)>v(\cdot;\sigma)$.
\begin{theorem}\label{thm:derive_zero}
Let Assumptions \textbf{A1-A3} for $\varsigma=-\alpha\ell$, where $\ell$ is given in \eqref{eqn:linear_coefficients} and $\alpha:[0,T]\times\mathbb{R}^2\to\mathbb{R}$ is a positive function. 
   If  $\nabla_\varsigma v(\cdot;\sigma)\equiv0$, then $v(\cdot,;\sigma)$ is the unique viscosity solution of  the fully nonlinear equation \eqref{eqn:fully_nonlinear} in the classs of functions with $p^\textrm{th}$ growth. 
\end{theorem}
\begin{eg}
    For \eqref{nonlinear_merton}, the semi linear equation is given by 
\end{eg}

\subsection{Gradient ascent based on the directional maximum principle}
\label{sec:gradient_ascent}
\begin{figure}
    \centering
    \includegraphics[width=0.8\textwidth]{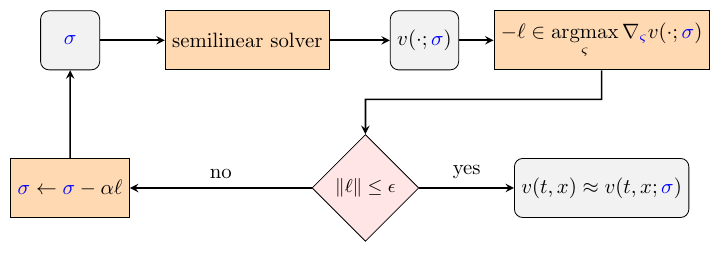}
    \caption{The flowchart of the scheme for fully nonlinear PDEs}
    \label{fig:max_principle}
\end{figure}

We start the scheme with initialization $\sigma(t,x):=\sigma_0(t,x)$. At step $m$ of the scheme, we solve the semilinear equation \eqref{eqn:semi_nonlinear} to obtain $v(\cdot,\sigma_m)$, from which we evaluate  $v(\cdot,\sigma_m)$, $\nabla v(\cdot,\sigma_m)$, and $D^2v(\cdot,\sigma_m)$, and consequently the coefficient $k_m$, $\mu_m$, and $\ell_m$ of the linear PDE \eqref{eqn: directional derivative of u}, given by \eqref{eqn:linear_coefficients}. Note that Theorem~\ref{thm:positive_gradient}  guarantees that ${\nabla}_{\varsigma} v(\cdot;\sigma_m)$ is  positive by  for  $\varsigma:=-\ell_m$.  
Then, we update $\sigma_{m+1}=\sigma_m-\ell_m$. We continue the scheme until 
$\sup_{t,x}{\nabla}_{\varsigma_n} v(t,x;\sigma_m)$ is sufficiently small. 
This indicates that the maximizer $\sigma^*$ is practically approximated and $v(\cdot;\sigma_m)$ is a good approximation of $v$.
\begin{remark}
    One should avoid solving the linear equation \eqref{eqn: directional derivative of u} to obtain the directional gradient. Instead, one can rely on the magnitude of $\ell_m$ to decide whether optimal $\sigma$ is attained or not. If $\sup_{t,x}\ell_m(t,x)<\epsilon$ for $\epsilon>0$ sufficiently small, then $\sup_{t,x}{\nabla}_{\varsigma_n} v(t,x;\sigma_m)<\epsilon T$. Therefore, the stopping  criterion of the algorithm can be solely based on $\ell_m$. 
    In addition, the value f ${\nabla}_{\varsigma_n} v(t,x;\sigma_m)$ is related to the value of $\sigma_m(s,x)$ on $s\in[t,T]$, which can acculumate small errors and does not necessarily represent how close $\sigma_m(t,x)$ is to the optimal $\sigma$ at the specific point $(t,x)$.
\end{remark}
This maximum principle scheme is summarized in Algorithm~\ref{alg:max_principle} and Figure~\ref{fig:max_principle}.

\begin{algorithm}[H]
    \caption{Maximum principle scheme}
    \label{alg:max_principle}
\begin{algorithmic}
\State{\textbf{Input}:
Fix $\epsilon>0$,  $M\in\mathbb{N}$, $\sigma_0:[0,T]\times Q\to\mathbb{M}_d$ , and a decreasing positive sequence $\{\alpha_n\}_{n\in\mathbb{N}}$.} 
\State{\textbf{Initialize}: Set $m=0$, $\sigma=\sigma_0$, and $\textsc{Gr}_0>\epsilon$. }
\While{$\textsc{B}_m > \epsilon$ \textit{and} $m \leq M$}

    \State{Solve semilinear equation \eqref{eqn:bsde_L-F} with $\sigma=\sigma_m$ using empirical risk minimization \eqref{eqn:risk_minimization_semilinear} to obtain approximation $\hat{v}(t,x;\sigma_m)$, $\nabla \hat{v}(t,x;\sigma_m)$. }

    \State{Compute $\ell_m$ in \eqref{eqn:linear_coefficients} as a function of $v=\hat{v}(t,x;\sigma_m)$ and its derivatives.}

    \State{$\sigma_{m+1}=\sigma_m-\alpha_m\varsigma_m$}

    \State{$m += 1$ and  $\textsc{B}_m=\|\ell_m\|$}
\EndWhile
\\
\Return{$\hat{\sigma}:=\sigma_m$ and $\hat{v}(\cdot;\hat{\sigma}):[0,T]\times Q\to\mathbb{R}$}
\end{algorithmic}
\end{algorithm}

\begin{remark}[Testing framework]
    In the light of \eqref{eqn:double_conjugate},
\eqref{eqn:minimizing_K}, and Theorem\ref{thm:derive_zero}, if $v(t,x)$ satisfies \eqref{eqn:fully_nonlinear}, for all $(t,x)$, we must have
\begin{equation}
\begin{split}
     \mathcal{E}[v](t,x):=&h\big(t,x,v(t,x),\nabla v(t,x),D^2v(t,x)\big)\\
     &-\frac12(\sigma^{\intercal}\sigma): D^2 v(t,x)+f\big(t,x,v(t,x),\nabla v(t,x),\sigma(t,x)\big)=0.
\end{split}
\end{equation}
Ideally, an approximate solution, $\hat{v}$ satisfies 
\begin{equation}
    \sup_{(t,x)}\mathcal{E}[\hat{v}](t,x)\le \epsilon
\end{equation}
for some small $\epsilon>0$. However, since our numerical scheme is based on the empirical risk, we only expect to have 
\begin{equation}
    \mathbb{E}\bigg[\int_0^T\mathcal{E}[\hat{v}](t,X_t)\bigg]\le \epsilon
\end{equation}
This can be used to verify that the approximate solution is well generalized along the sample paths that are not in the training data and detect overfitting. In the current paper, we used the closed form solution to verify that the trained solution is not overfitting to the data.
\end{remark}

\subsection{Numerical solution for the semilinear equation} 
\label{sec:numerical semilinear}
In Algorithm~\ref{alg:max_principle}, each iteration involves solving semilinear problem \eqref{eqn:num_semilinear}. Theoretically, any numerical method for the semilinear problem can be used as long as it provides us with an approximately accurate estimation of $\ell$ in \eqref{eqn:linear_coefficients}. However, some methods, such as finite difference, may be sensitive to dimension and cost intensive. For the gradient ascent method, our preferred method is a modification of the one in \citet{JHW17}, which provides the solution to the semilinear equation via Monte Carlo simulation of $dX_t=\sigma(t,x_t)dB_t$ and an empirical risk minimization for deep neural networks. On of the advantages of this method is the use of automatic differentiation to find the derivatives of the solution efficiently. In addition, a neural network provides a global estimate for solution and its derivatives, which saves us an interpolation step. Since the method used for semilinear problem is not pivotal to the main result of this paper, we describe it the in Appendix~\ref{sec:semi_linear_scheme}. 
\subsection{Implementation consideration}
Based on our preferred method discussed in Section~\ref{sec:numerical semilinear}, the approximation of $\ell_m$ in Algorithm~\ref{alg:max_principle} is a random variable. As a matter of fact, the Monte Carlo samples of $dX_t=\sigma(t,x_t)dB_t$ can lie outside the region where the approximation of solution and its derivatives with neural network is accurate. Therefore, the approximation of $\ell_m$ can possible take large values. To add stability to the gradient ascent method, we introduce a bound, $b_m$ as  hyperparameter and replace $\ell_m$ in Algorithm~\ref{alg:max_principle} with $-b_m\vee(\ell_m\wedge b_m)$. For the same reason, an $L^p$-norm for  $\textsc{B}_m$ is the stopping criteria is more appropriate, $\textsc{B}_m=\|\ell_m\|_p$. In our implementations, we considered $\textsc{B}_m=\|\ell_m\|_1$. To increase the stability of the method, we also eliminated the outliers from the simulated paths of $dX_t=\sigma(t,x_t)dB_t$. The approximation of the value function and its derivatives at thee outliers grossly inaccurate and lead to wrong overall approximation.

\section{Numerical Experiment}\label{sec:numerics}
The portfolio selection problem with stochastic volatility, introduced in \citet{Chesney-Scott-1989}, is chosen to host the numerical experiment. This specific choice is also made to compare our results to \citet{GPW21}, where they used a method based on step-wise dynamic programming to solve a fully nonlinear problem.

Consider $d$ risky assets with price process $S=(S^1,...,S^d)$ given by 
\begin{equation}
    dS_{t}=\textrm{\normalfont diag}(S_{t})\textrm{\normalfont diag}(\exp(Y_t))(\lambda(Y_t) dt+dB^{0}_{t}) 
\end{equation}
where $\textrm{\normalfont diag}(a)$ is the diagonal matrix with elements of vector $a$ on the diagonal, $B^0_t$ is a $d$-dimensional Brownian motion, $\lambda:\mathbb{R}^d\to\mathbb{R}^d$ is the risk premium of the assets with $\lambda(y) = (\lambda_i(y_i))^\intercal$, exponential is element-wise, and $Y_t$ is an $\mathbb{R}^d$-values OU modeling the stochastic volatility:
\begin{equation}\label{eqn:OU}
    dY_{t}=\textrm{\normalfont diag}(\kappa)(\theta-Y_t)dt+\textrm{\normalfont diag}({\nu}) dB^{1}_{t}
\end{equation}
where $\kappa$, $\theta\in\mathbb{R}^d$, $\sigma\in\mathbb{R}_+^d$, and $B_t$ is a $d$-dimensional Brownian motion such that $d\langle B^0,B^1\rangle =\textrm{\normalfont diag}(\rho) dt$, with $\rho\in(-1,1)^d$.   Consider an admissible trading strategy,  $\pi=\{\pi_t\}_{t\ge0}$; where $\pi_i \in \mathbb{R}^d$ is the amount invested in  asset $i$ at time $t$. In this model, the wealth process is given by
\begin{equation}
    d W^\pi_t =\pi_{t}^\intercal(\textrm{\normalfont diag}(\exp(Y_t))(\lambda(Y_t) dt+dB^0_{t})
\end{equation}
The objective of the agent is to maximize his expected utility from terminal wealth over all admissible trading strategies $\pi$. Hence, the value function $u$ is given by: 
\begin{equation}
    u(t,x)=\sup_{\pi \in \mathcal{A}} \mathbb{E}[U(W^\pi_T)]
\end{equation}
where $x=(w,y_1,...,y_d)$ and $\mathcal{A}$ is the collection of all trading strategies such that $W_t^\pi\ge0$ for all $t\in[0,T]$ a.s.
The Hamilton-Jacobi-Bellman (HJB) for value function $u$ the above is given by
\begin{equation}\label{eqn:HJB}
\left\{\begin{matrix}
\begin{aligned}
& 0 =-\partial_tv - \mathcal{L}^y v- h\big(w,y,\nabla v(t,w,y),D^2 v(t,w,y)\big) \\
\\ 
&v(T,x)=U(w)
\end{aligned}
\end{matrix} \right.
\end{equation}
where 
\begin{equation}
    h(w,y,q,\gamma):=-
\dfrac12\dfrac{\textstyle\sum_{i=1}^{d}\big(\lambda_i(y_i)q_0+\rho_i\nu_i\gamma_{0i}\big)^2}{\gamma_{00}}
\end{equation}
with $\mathcal{D}_h=\{\gamma:\gamma_{00}<0\}$ and $\mathcal{L}^{v}=\textstyle\sum_{i=1}^{d}[\kappa_i(\theta_i-y_i)\partial_{y_i}+\frac{1}{2}{\nu}_i^2\partial_{y_i}^2]$ is the infinitesimal generator of the \eqref{eqn:OU}. 
The Legendre transform of $h$ is given by
\begin{equation}\label{eqn: F when n=d}
    f(t,x,p,q,\sigma)=
        -|\sigma_{00}||q_0|\sqrt{\textstyle\sum_{i:\rho_i=0}|\lambda_i(y_i)|^2}-|\sigma_{00}|\sigma_{01}q_0\textstyle\sum_{i:\rho_i\neq0}\frac{\lambda_i(y_i)}{\rho_i\nu_i}
\end{equation}
The domain of $f$, $\mathcal{D}_f$ consists of all lower triangular $d+1$ by $d+1$ matrices $\sigma=[\sigma_{ij}:0\le i\le j\le d]$ such that 
\begin{equation}
    |\sigma_{0i}|\le \rho_i\nu_i, ~\sigma^2_{0i}+\sigma^2_{ii}=\nu_i^2,~\textstyle \textstyle\sum_{0\le k\le i\wedge j}\sigma_{i,k}\sigma_{j,k}=0, \textrm{ and }~\textstyle\sum_{0\le k\le i}\sigma_{i,k}^2=\nu_i^2
\end{equation}
We simplify the problem by taking $\lambda_i(y_i) = \lambda_i y_i$ for $\lambda_i\in\mathbb{R}$ and the utility function $U(x)=1-e^{-\eta x}$, for which a closed-form solution is given by 
\begin{equation}\label{eqn:closed-formSC}
    v(t,x,y)=1-\exp\Big({-\eta w-\textstyle\textstyle \textstyle\sum_{i=1}^d\phi_i(T-t)y_i^2-\psi_i(T-t)y_i-\chi_i(T-t)}\Big)
\end{equation}
where $\phi_i$, $\psi_i$, and $\chi_i$ are defined in \citet[Appendix]{SZ99} or \citet[Section 3.5 page 16]{GPW21}.

In the no-leverage case,  $\rho=0$, we have  
\begin{equation}\label{eqn: H}
        h(q,\gamma) =\tfrac{q^2_0}{2\gamma_{00}}\Lambda^2(y),~~~\textrm{ with }~~~\Lambda(y)=\big(\textstyle\sum_{i=1}^d\lambda_i^2(y_i)\big)^{\sfrac12}
\end{equation}
and Legendre transform of $h$ is given by
\begin{equation}\label{eqn: F when n=d and rho = 0}
f(t,x,p,q,\sigma)=-\textstyle\sum_{i=1}^{d}\left[\kappa_i(\theta_i-y_i)q_i\right]- |q_0||\sigma_{00}|\Lambda(y)\\
\end{equation}
With domain of $f$ is made of matrices of the form $\textrm{diag}(\sigma_{00}(t,w,y),\nu_1,...,\nu_d)$ and therefore, in only varies in $\sigma_{00}$.
Therefore, the maximizer $\sigma^*$ for \eqref{eqn:double_conjugate} is given by  
\begin{equation}
    \sigma_{00}^*(t,w,y)=\mathop\textrm{\normalfont argmax}\limits_{\sigma\in\mathcal{D}_f}\Big\{\frac12\sigma^2_{00}\partial_{ww}v+\Lambda(y)|\partial_{w}v ||\sigma_{00}|\Big\}
\end{equation}
By using \eqref{eqn:closed-formSC}, we obtain $\sigma_{00}^*(t,y)=\tfrac{\Lambda(y)}{\eta}$ and  
\begin{equation}\label{eqn: optimal quadratic}
    \sigma^*=\textrm{\normalfont diag}\left(\sigma_{00}^*,\nu_1\cdots\nu_d\right)
\end{equation}

\subsection{Reduction to Merton problem}
For the choice of $\rho_i=0$ for all $i$ and constant $\lambda_i(y_i) = \lambda_i\in\mathbb{R}$, the problem reduces to Merton's problem and the closed form is given by $v(t,x)=1-\exp\big({-\eta w-\frac12\Lambda(T-t)}\big)$ with $\Lambda=\textstyle\textstyle \big(\textstyle\sum_{i=1}^d\lambda_i^2\big)^{\sfrac12}$.

We implement the gradient ascent on this reduced problem  with $\sigma_0=\textrm{\normalfont diag}\left(\sigma_0,\nu_1\cdots\nu_d\right)$. For  a random choice of $\sigma_0(t,x)\equiv\sigma_0\in\mathbb{R}_+$. Given $\sigma_n(t,x)$, at  iteration $n$, we solve the semilinear PDE
\begin{equation}\label{eqn:num_semilinear}
\begin{cases}
      0=-\partial_t v_n-\frac12\sigma^2_n \partial_{ww} v_n -\mathcal{L}^y v_n -\Lambda(y) |\partial_w v_n|\sigma_n\\
        v_n(T,\cdot)=U(x)
\end{cases}
\end{equation}
where $\Lambda(y):=\big(\textstyle\sum_{i=1}^d\lambda_i^2(y_i)\big)^{\sfrac12}$.
In \eqref{eqn:linear_coefficients}, we evaluate
\begin{equation}
    \begin{split}
        k(t,x)&\equiv0\\
        \mu_n(t,x)&=\Big(|\lambda|\textrm{sgn}({\partial_w v_n}(t,x))\sigma_n(t,x),\kappa_1(\theta_1-y_1),...,\kappa_d(\theta_d-y_d)\Big)\\
        \ell_n(t,x)&=[L_{ij}(t,x):i,j=0,...,d]
    \end{split}
\end{equation}
where $L_{ij}(t,x) \equiv 0$ except for
\begin{equation}
    L_{00}(t,x) = -\Lambda(y) |\partial_w v_n(t,x)|-\sigma_n(t,x) \partial_{ww}v_n(t,x)
\end{equation}
 
Then, we set $\varsigma_n=-\alpha_n\ell_n$, for some $\alpha_n\in\mathbb{R}_+$ and update $\sigma_{n+1}(t,x)=\sigma_n(t,x)-\alpha_n\ell_n(t,x)$.

The detained implementation can be found at \citet{deeprepo}.

\subsection{Neural network architecture}
  The neural network for the solver of semilinear problem consist of one input layer, three hidden layers, and one output layer. The activation function utilized is $\tanh$. The number of neurons in each hidden layer has been chosen dependent of the dimension of the problem. Since we incorporated the time variable $t$ into the input, the input layer has  dimension $d+2$, $(t,w,y_1,...,y_d)$.  The neural networks employed to approximate $v(\cdot;\sigma)$ have an output dimension of $1$.

For semilinear equation \eqref{eqn:num_semilinear}, the stochastic process is given by $X=(Y^0,...$ $,Y^d)$, where  $dY^0_t=\sigma_n(t,X)dB^0_t$ and,  for  $Y^i$, $Y^i$ are given by \eqref{eqn:OU}. The stochastic process is given by $X=(Y^0,..,Y^d)$, where  $dY^0_t=\sigma_n(t,X)(|\lambda|\textrm{sgn}({\partial_w v_n}(t,x))$ $dt+dB^0_t)$ and, for  $Y^i$, are given by \eqref{eqn:OU}. After generating sample paths, the risk function in \eqref{eqn:risk_minimization_semilinear} will be replaced by an empirical risk function. The number of sample paths for each problem is recorded in Table~\eqref{tab:universal dimnesion 3, 5 and 10} along with the choice of minimizer and the number of epochs.

\subsection{Single asset \texorpdfstring{$d=1$}{d=1}}
\begin{figure}[htbp] 
    \centering
    \subfigure[]{\includegraphics[width=0.28\textwidth]{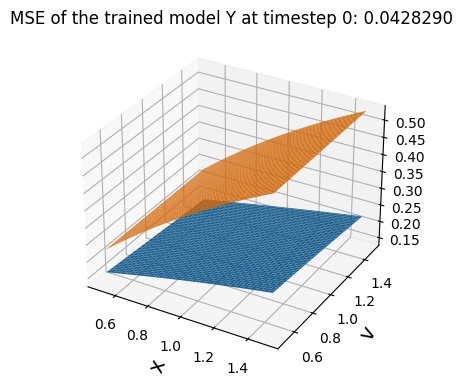}} 
    \subfigure[]{\includegraphics[width=0.28\textwidth]{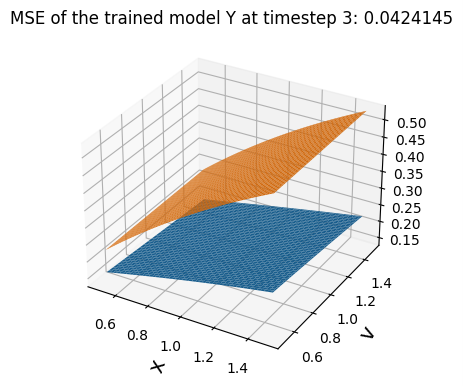}} 
    \subfigure[]{\includegraphics[width=0.28\textwidth]{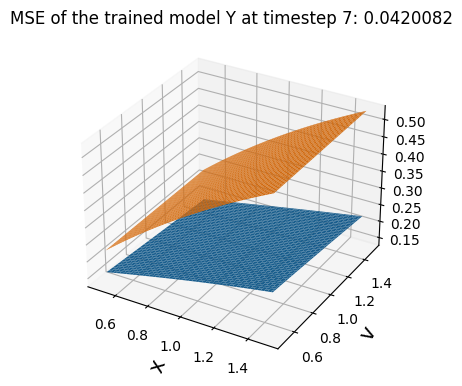}}
    \caption{$d=1$. The value $v_0(t,x,y)$ for a randomly initialized $\sigma_0$ at $t=0, t_3$, and $t_7$.}
\label{fig:Initial}
\end{figure}

\begin{figure}[htbp] 
    \centering
    \subfigure[]{\includegraphics[width=0.28\textwidth]{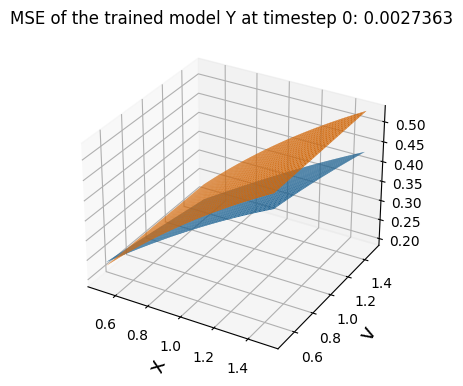}} 
    \subfigure[]{\includegraphics[width=0.28\textwidth]{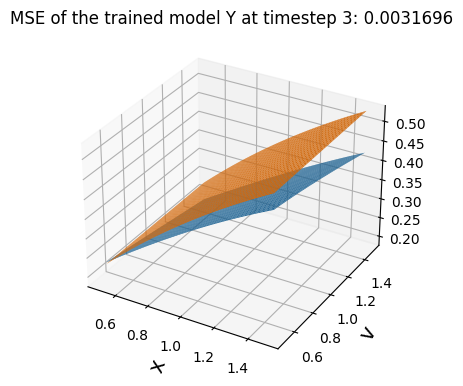}} 
    \subfigure[]{\includegraphics[width=0.28\textwidth]{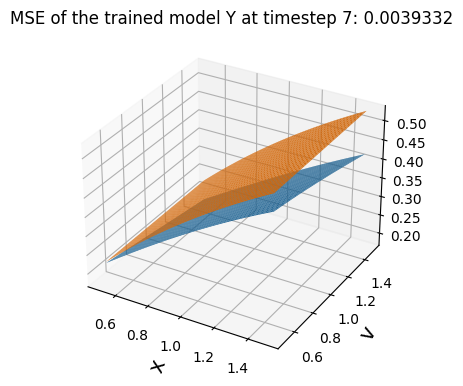}}
    \caption{The solution at time steps 0, 3, and, 7}
\label{fig:5th update}
\end{figure}

\begin{figure}[htbp] 
    \centering
    \subfigure[]{\includegraphics[width=0.28\textwidth]{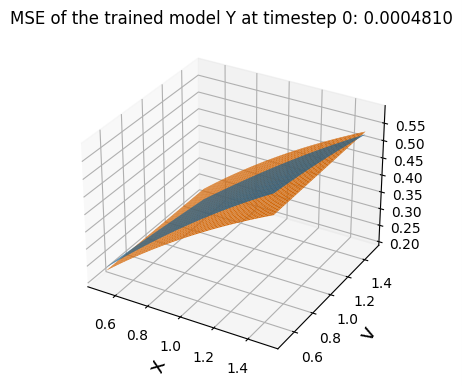}} 
    \subfigure[]{\includegraphics[width=0.28\textwidth]{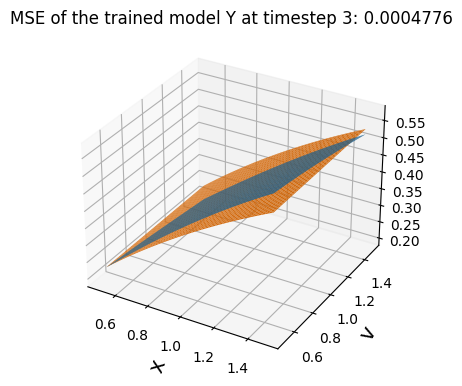}} 
    \subfigure[]{\includegraphics[width=0.28\textwidth]{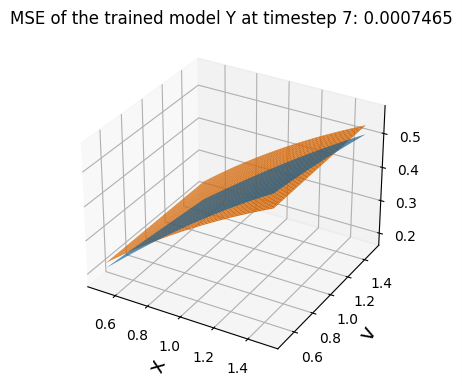}}
    \caption{The solution $v_{10}$ at time steps 0, 3, and, 7}
\label{fig:10th update}
\end{figure}

When there is only one asset with stochastic volatility, the problem has dimension two. We choose $T=1$ and discretize time by $\Delta t = \frac{T}{20}$ and $t_n=n\Delta t$. The results of scheme for the initial choice of $\sigma_0$, $\sigma_3$, and $\sigma_7$ for different time steps are are illustrated in Figures \eqref{fig:Initial}, \eqref{fig:5th update}, and \eqref{fig:10th update}. These figures serve to demonstrate the scheme's efficacy by showcasing the convergence of the approximated solution (depicted in blue) towards the true solution (depicted in orange) as the number of iterations (sigma updates) increases across different time steps. We conducted our testing using a sample size of $2^{13}$ and a batch
size of $2^8$. 
The figures \ref{fig:Initial}-\ref{fig:10th update} show that the approximated solution at different time steps approach the true solution as the number of iteration (updates of $\sigma_m$) increases. At the $10^{th}$ iteration, \ref{fig:10th update}, the error at randomly selected points at $t=0, t_3$, $t_7$ is measured to be 0.0002838, 0.0005062, and 0.0003246, respectively.  The error is evaluated compared to the closed-form solution \eqref{eqn:closed-formSC}.

\subsection{Multiple assets \texorpdfstring{$d>1$}{d>1}}
Table \ref{tab:universal dimnesion 3, 5 and 10} provides a summary of the hyperparameters and the results obtained for higher dimensions when $d=2,4,9$. We denote $\theta$ and $\beta$ the set of parameters for training the neural networks that approximate solutions of equation \eqref{eqn:num_semilinear}, respectively.  The number of time steps used for each dimension is the same as $d=1$, $20$. The evaluation metrics used in the analysis of higher dimension are the mean square error (MSE) and the standard deviation (STD). The MSE and STD values were calculated based on $10^3$ randomly sampled points within the designated testing regions. The reported results are averaged over three independent runs, with each run consisting of 15 iterations. The learning rate used for each optimizer is  $10^{-3}$ and the number of epochs are $50$.

\section{Conclusion}
In this study, we showed that the idea of reducing a fully nonlinear problem to a class of semilinear problems, introduced in \citet{STZ12}, can be used as a basis for a numerical approximation. In particular, we showed plausibility of transferring efficient methods for semilinear problems to relatively efficient methods for fully nonlinear problems. While we established the convergence, the accuracy depends on the accuracy of the scheme for semilinear equations. 

As discussed in \citet{HD23}, the time of the algorithm proposed in this study in not sufficiently fast, due to slow solvers for semilinear PDEs. Increase in time-efficiency of numerical schemes for semilinear PDEs is the subject of a different line of research in our future endeavor. 

\vskip 5pt  \begin{table}[htbp!]
\centering
\caption{Summary of Hyperparameters Used in Higher Dimensions}
\begin{tabular}{|>{\centering}m{0.2cm}|c|c|>{\centering}m{0.8 cm}>{\centering}m{0.8 cm}|>{\centering}m{2 cm}|c|}
\hline
\multirow{2}{*}{$d$}& \multirow{2}{*}{Sample size}&\multirow{2}{*}{Batch size} &\multicolumn{2}{c|}{Optimizers} &\multirow{2}{*}{Avg MSE} & \multirow{2}{*}{Avg STD} \\ 

& & & $\theta$ & $\beta$ & &  \\ \hline\hline

1         & $2^{13}$ &   $2^{8}$               & Adam   & SGD  & $3.411\times10^{-3}$     & $3.21\times10^{-3}$       \\ \hline
2         & $2^{14}$ &   $2^{11}$               & Adam   & SGD  & $3.411\times10^{-3}$     & $3.21\times10^{-3}$       \\ \hline
4        & $2^{16}$  &   $2^{13}$                &SGD   & SGD    & $2.61\times10^{-3}$   &$2.1\times10^{-4}$   \\ \hline
9       & $2^{18}$  &   $2^{15}$                  &SGD  & SGD          & $3\times10^{-4}$   & $3.3\times10^{-3}$          \\ \hline
\end{tabular}
\label{tab:universal dimnesion 3, 5 and 10}
\end{table}

\begin{figure}[htbp] 
    \centering
    {\includegraphics[width=0.32\textwidth,height=0.15\textheight]{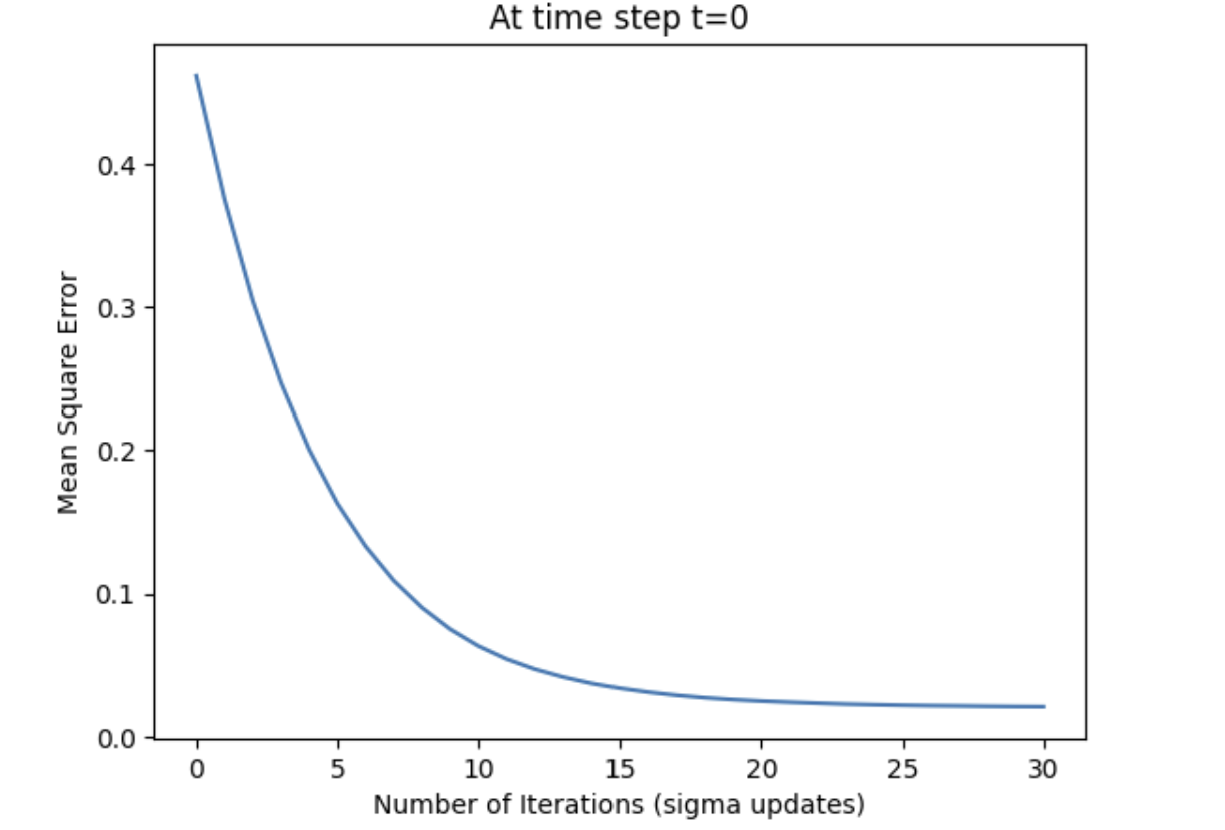}} 
    {\includegraphics[width=0.32\textwidth,height=0.15\textheight]{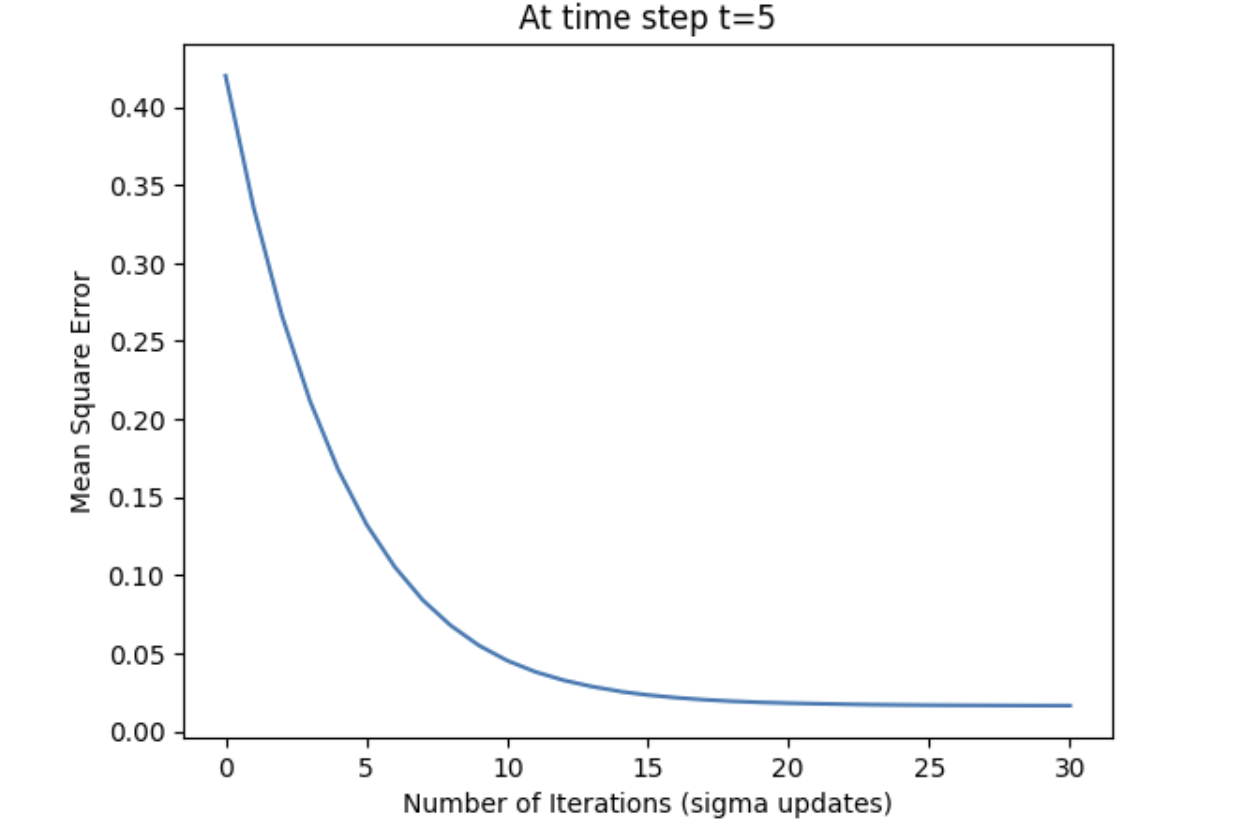}} 
    {\includegraphics[width=0.32\textwidth,height=0.15\textheight]{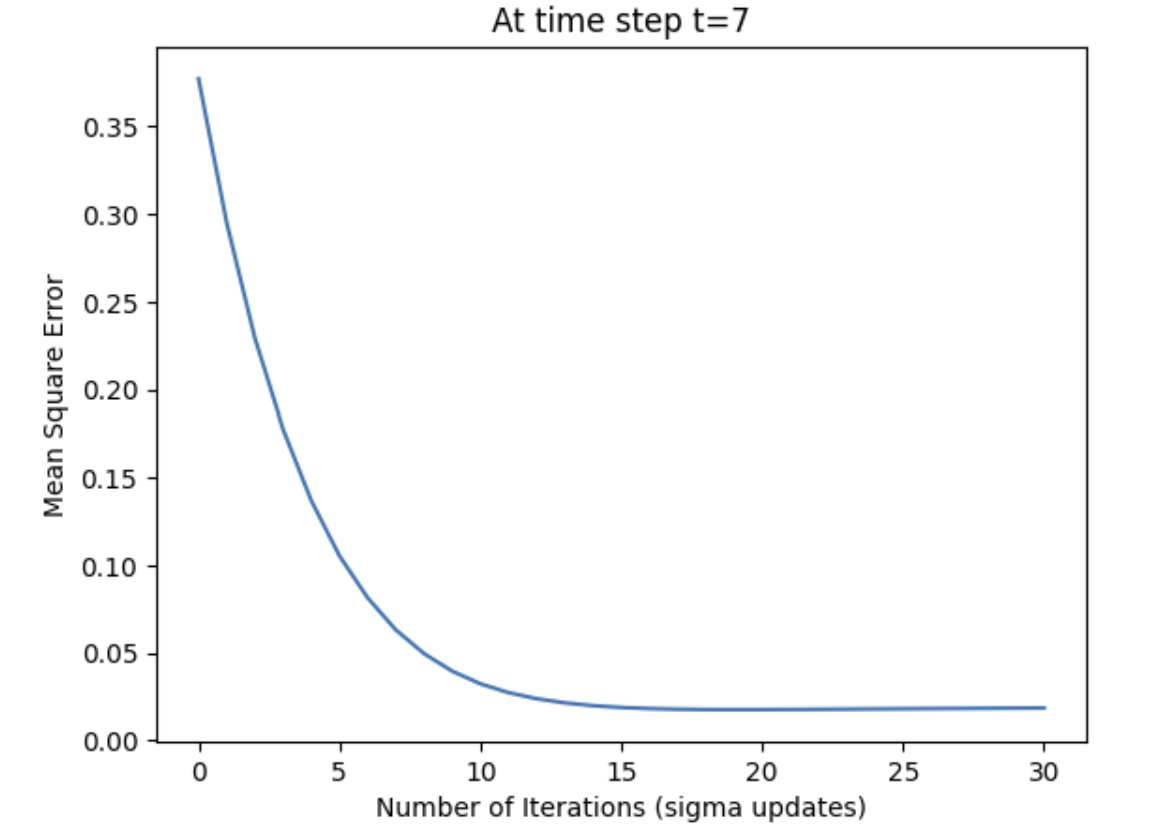}}
    \caption{MSE vs $M$, number of iterations for $\sigma$, at $t_0=0. t_5$, and $t_7$ for $d=2$}
\label{fig:universal dimnesion 3}
\end{figure}

\begin{figure}[htbp] 
    \centering
    {\includegraphics[width=0.32\textwidth,height=0.15\textheight]{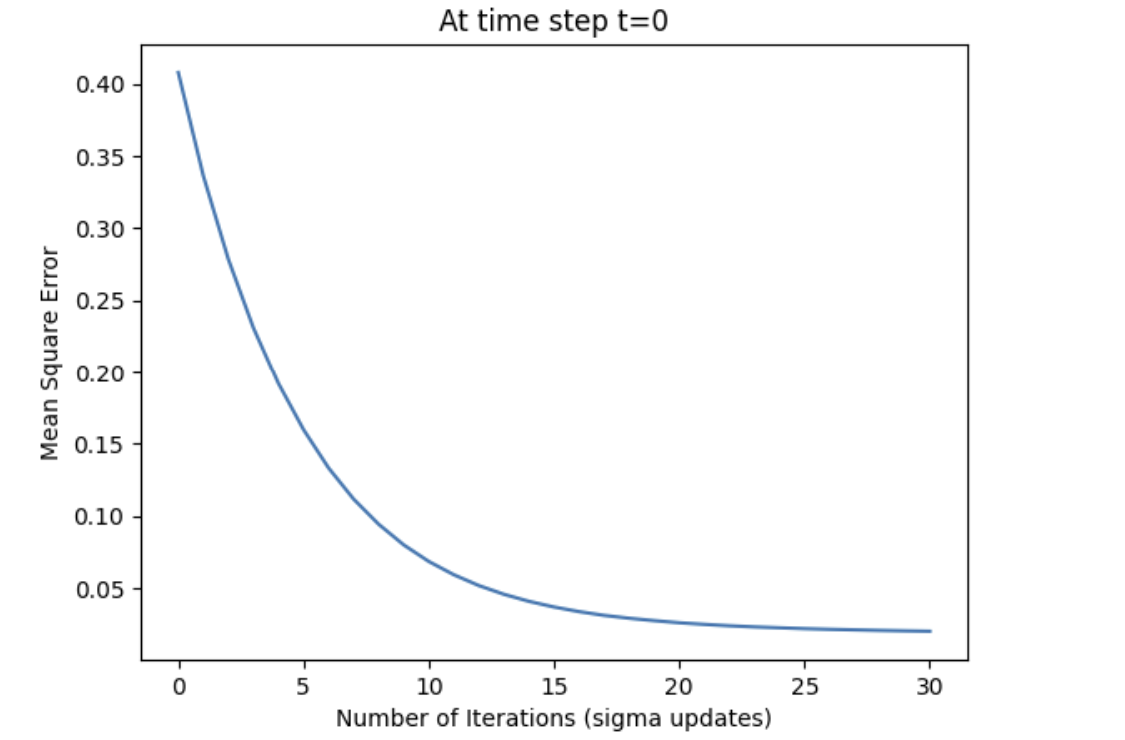}} 
    {\includegraphics[width=0.32\textwidth,height=0.15\textheight]{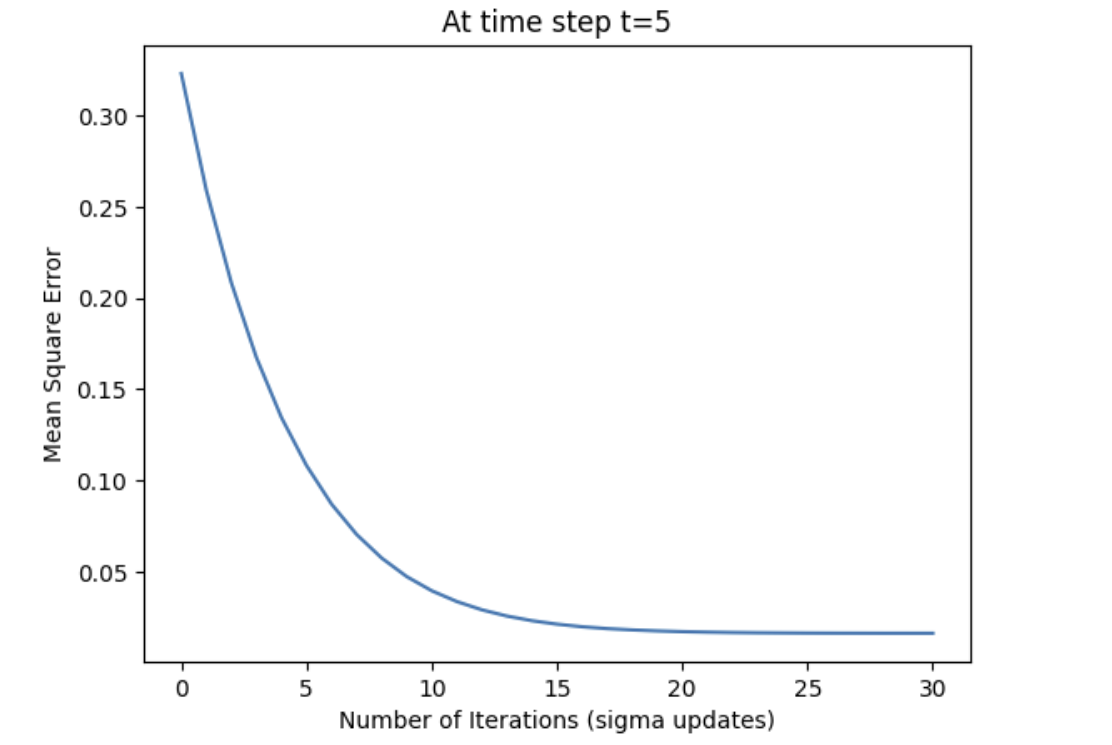}} 
    {\includegraphics[width=0.32\textwidth,height=0.15\textheight]{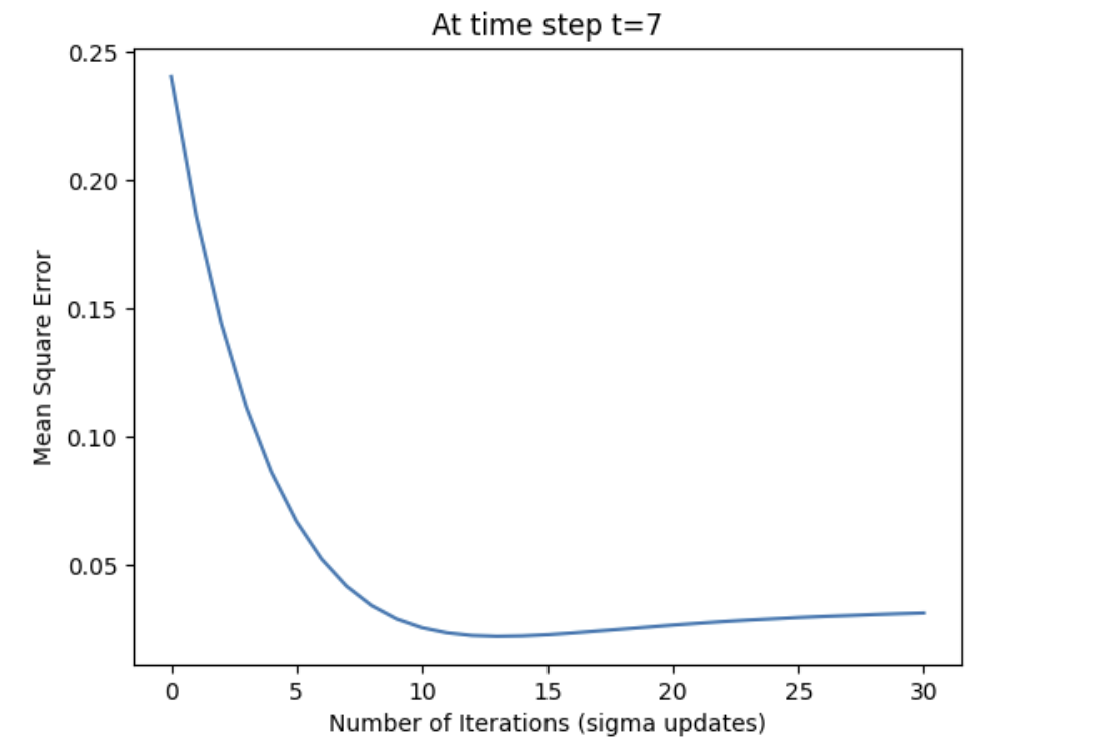}}
    \caption{MSE vs $M$, number of iterations for $\sigma$, at $t_0=0. t_5$, and $t_7$ for $d=4$}
\label{fig:universal dimnesion 5}
\end{figure}

\begin{figure}[htbp] 
    \centering
    {\includegraphics[width=0.32\textwidth,height=0.15\textheight]{figures/universal_MSE_t0.png}}
    {\includegraphics[width=0.32\textwidth,height=0.15\textheight]{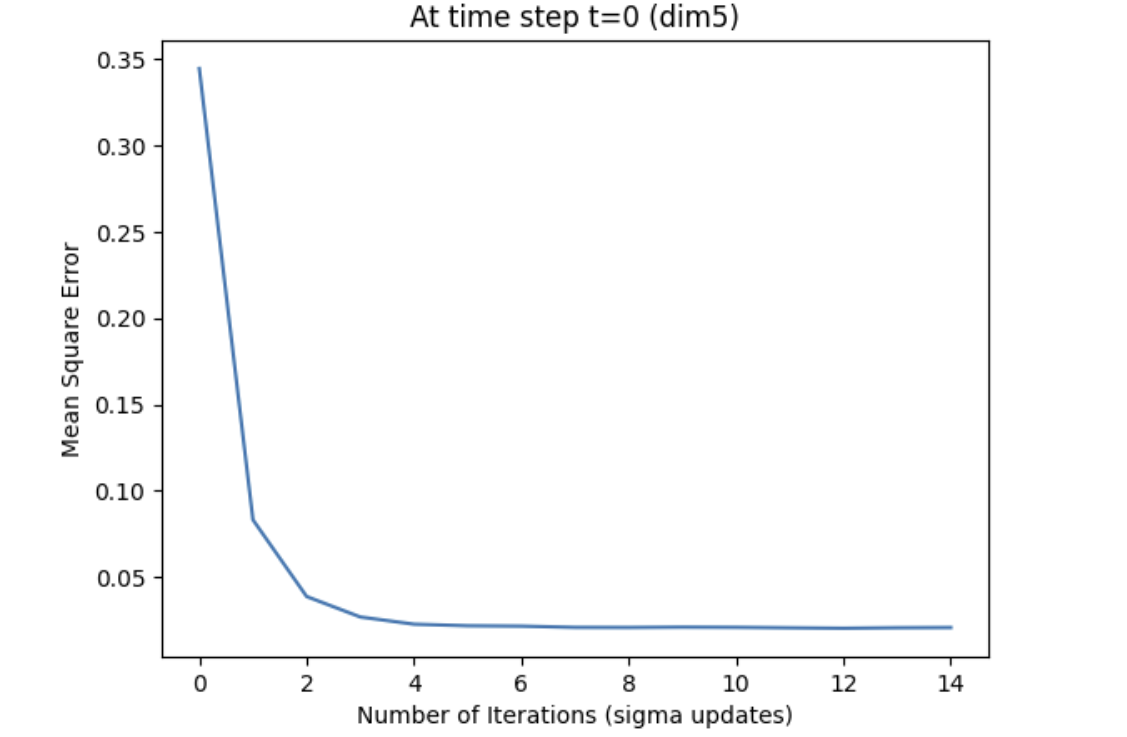}} 
{\includegraphics[width=0.32\textwidth,height=0.15\textheight]{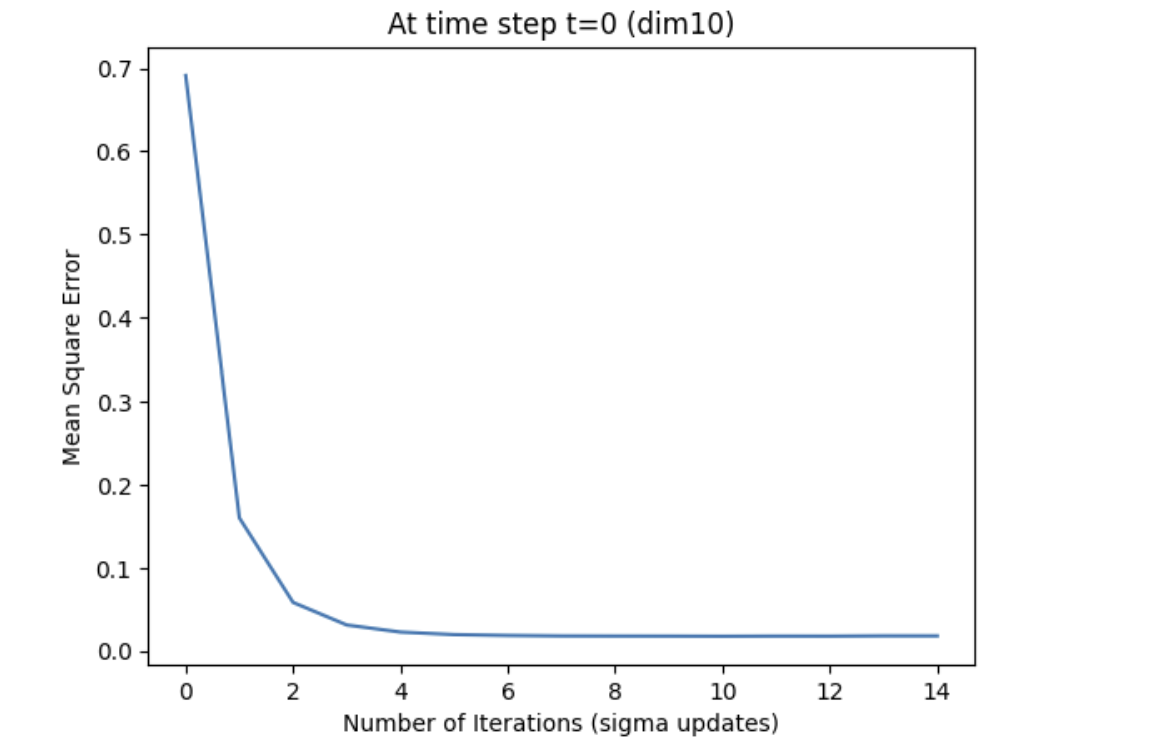}} 
    \caption{MSE vs $M$, number of iterations for $\sigma$, at $t_0=0$ for $d=2, 4$, and $9$}
\label{fig:universal dimnesion 5 and 10}
\end{figure}

\newpage
\appendix 
\section{Proof of the main results}
\label{sec:appendix}
The main results of this paper relies on the comparison principle for the parabolic problems, i.e., \eqref{eqn:fully_nonlinear}, \eqref{eqn:semi_nonlinear}, or \eqref{eqn: directional derivative of u}. For such nonlinear or linear PDEs, the general comparison principle is often established in the viscosity sense, a weak form of solutions. For the sake of completeness, we provide the basic definitions in this appendix and refer the reader \citet{CIL92} or \citet{T12} for more details on the theory of viscosity solutions.  

Consider the following problem.
 \begin{equation}\label{eqn:general}
 \begin{cases}
 	 	-\partial_t{u}-G(t,x,u(t,x),\nabla u(t,x),D^2u(t,x),\sigma(t,x))=0& (t,x)\in Q\\
 	 	u(T,x)=g(x)&x\in \mathbb{R}^d
 \end{cases}
 \end{equation} 
 Here, $G(t,x,p,q,\gamma,\sigma)$ is given by $h(t,x,p,q,\gamma)$ for problem \eqref{eqn:fully_nonlinear} and by $\frac12(\sigma^\intercal\sigma): \gamma-f(t,x,p,q,\sigma)$ for problem \eqref{eqn:semi_nonlinear}. Here, we assume $G$ is a continuous function.
The definition of viscosity solutions is based on replacing the derivative in the PDE with derivative of smooth functions tangent to the solutions at each point, as defined below. 
 \begin{defn}\label{def:viscosity}
    A upper semi-continuous (resp. lower semi-continuous) function $v$ (resp. $V$) a viscosity sub- (resp. super-)solution of  \eqref{eqn:general}  if 
    \begin{enumerate}[label=\arabic*)]
        \item if for any $(t,x)\in Q$ and for any smooth function $\phi:[0,T)\times\mathbb{R}^d\to\mathbb{R}$ such that $0=\max_{s,y}(v-\phi)(s,y)=v(t,x)-\phi(t,x)$ (resp. $0=\min_{s,y}(V-\phi)(s,y)=V(t,x)-\phi(t,x)$)
        \begin{equation}
            -\partial_t{\phi}(t,x)-G(t,x,\phi(t,x),\nabla \phi(t,x),D^2\phi(t,x),\sigma(t,x))\le 0 \textrm{(resp. $\ge 0 $)}
        \end{equation}
        \item $v(T,x)\le g(x)$ (resp. $V(T,x)\ge g(x)$)
    \end{enumerate}
  A (continuous) function which is both sub and super solution is called a viscosity solution.
\end{defn}
\begin{remark}
    In the definition above, the test function $\phi$ only matters through values $\partial_t{\phi}(t,x)$, $\phi(t,x)$, $\nabla \phi(t,x)$ ,$D^2\phi(t,x)$. Therefore, modifying $\phi$ such that this values are preserved does not affect the definition of viscosity solution. In particular, one can replace $\phi(s,y)$ by $\phi(s,y)+\frac{\lambda^2}{2}|y-x|^2$ to assures that $(t,x)$ is the unique global maximizer of $v-\phi$. 
\end{remark}
We also require the comparison principle for viscosity solutions to establish our results. More precisely, Definition~\ref{def:regularity} and the heuristic argument in Section~\ref{sec:intro} rely on the comparison principle for viscosity solutions.
\begin{defn}
    \label{def:comparison}
    We say that problem \eqref{eqn:general} satisfies comparison principle in the class of functions with $p^\textrm{th}$ growth condition, if for any upper semi-continuous viscosity subsolution $v$ and any  lower semi-continuous viscosity supersolution $V$ with 
    $|v(t,x)|+|V(t,x)|\le C(1+|x|^p)$, we have 
    $v(t,x)\le V(t,x)$ for all $(t,x)\in[0,T]\times Q$.
\end{defn}
Let $v(\cdot;\sigma)$ be the solution to \eqref{eqn:semi_nonlinear}. For $\sigma\in\mathcal{D}_f$ and $\varsigma:[0,T]\times Q\to\mathbb{M}_d$, define 
\begin{equation}\label{eqn:u^*u_*}
    \begin{aligned}
    u^*(t,x):=\limsup_{(s,y,\epsilon)\to(t,x,0)}\frac{v(s,y;\sigma+\epsilon \varsigma)-v(s,y;\sigma)}{\epsilon}\\
    u_*(t,x):=\liminf_{(s,y,\epsilon)\to(t,x,0)}\frac{v(s,y;\sigma+\epsilon \varsigma)-v(s,y;\sigma)}{\epsilon}
\end{aligned}
\end{equation}
We start the proof with the following lemma.
\begin{lemma}\label{lem:u^*_p_growth}
    Let Assumptions \textbf{A2-A3} holds. Then, $u^*$ and $u_*$ in \eqref{eqn:u^*u_*} are lower semicontinuous and upper semi-continuous,  respectively. In addition, $|u^*(t,x)|+|u_*(t,x)|\le C(1+|x|^p)$.
\end{lemma}
\begin{proof}
    By Assumptions~\textbf{A2-A3} and Remark~\ref{rem:A2}, there exists $C>0$ and $p\ge0$ independent of $\epsilon\in[0,\epsilon_0]$ such that
    \[
    \dfrac{|v(s,y;\sigma+\epsilon \varsigma)-v(s,y;\sigma)|}{\epsilon}\le C(1+|y|^p)
    \]
    which readily shows the bound on $u^*$ and $u_*$. The semi-continuity of $u^*$ and $u_*$ is the direct result of the definition of $\limsup$ and $\liminf$.
\end{proof}
\begin{lemma}\label{lem:u^*_sub}
Let Assumptions \textbf{A2-A3} holds.    $u^*$ (resp. $u_*$) is a sub (resp. super) solution of \eqref{eqn: directional derivative of u}.
\end{lemma}
\begin{proof}
    We only show the sub solution property as the super solution property follows the same lines of argument. For any $(t,x)\in Q$ let  smooth function $\phi:[0,T)\times\mathbb{R}^d\to \mathbb{R}$ such that $0=\max_{s,y}(u^*-\phi)(s,y)=(u^*-\phi)(t,x)$ and $(t,x)$ is the unique global maximizer. As a property of $\limsup$, there exists a sequence $\{(s_n,y_n,\epsilon_n)\}_n$ with $(s_n,y_n,\epsilon_n)\to(t,x,0)$ such that $\mathop{\lim}\limits_{n\to\infty}\frac{v(s_n,y_n;\sigma+\epsilon_n \varsigma)-v(s_n,y_n;\sigma)}{\epsilon_n}=u^*(t,x)$. In particular, one can take $(s_n,y_n)$ to be maximizer of $v(s,y;\sigma+\epsilon_n \varsigma)$ $-v(s,y;\sigma)-{\epsilon_n}\phi(s,y)$. This is achieved by modifying the test function $\phi$ such that the maximizer of $u^*-\phi$ is unique and and the maximizer of $\dfrac{v(s,y;\sigma+\epsilon_n \varsigma)-v(s,y;\sigma)}{\epsilon_n}-\phi(s,y)$, $({t}_n,{x}_n)$, is attained inside a ball centered at $(t,x)$. Then, up to a subsequence, $(t_n,x_n)\to(\bar{t},\bar{x})$ and the following calculation shows that $(\bar{t},\bar{x})$ is a maximizer of $(u^*-\phi)$.
    \[
    \begin{aligned}
        0&=(u^*-\phi)(t,x)=(u^*-\phi)(\bar{t},\bar{x})\\
        &\ge \limsup_{n\to\infty}\dfrac{v(t_n,x_n;\sigma+\epsilon_n \varsigma)-v(t_n,x_n;\sigma)}{\epsilon_n}-\phi(t_n,x_n)\\
        &\ge \lim_{n\to\infty}\frac{v(s_n,y_n;\sigma+\epsilon_n \varsigma)-v(s_n,y_n;\sigma)}{\epsilon_n} -\phi(s_n,y_n)=(u^*-\phi)(t,x)=0
    \end{aligned}
    \]
    Since the maximizer is unique, $(t,x)=(\bar{t},\bar{x})$, $(t_n,x_n)\to({t},{x})$, and  
    \[\lim_{n\to\infty}\frac{v(t_n,x_n;\sigma+\epsilon_n \varsigma)-v(t_n,x_n;\sigma)}{\epsilon_n}=u^*(t,x)
    \] 
    Set $c_n:=\frac{v(t_n,x_n;\sigma+\epsilon_n \varsigma)-v(t_n,x_n;\sigma)}{\epsilon_n}$. By the above argument, $(t_n,x_n)$ is the maximized of $v(s,y;\sigma+\epsilon_n \varsigma)-\big(v(s,y;\sigma)+\epsilon_n\phi(s,y)+\epsilon_nc_n\big)$. In particular, since $\sigma$ is regular, in the sense of Definition~\ref{def:regularity}, $v(s,y;\sigma)+\epsilon_n\phi(s,y)+\epsilon_nc_n$ constitutes a test function for $v(s,y;\sigma+\epsilon_n \varsigma)$. Therefore, by assumption of continuous dependence of $\sigma$ in the direction of $\varsigma$, we have
    \begin{equation}
        \begin{aligned}
         0\ge &-\partial_{t}v(t_n,x_n;\sigma)-\frac{1}{2}(\sigma+\epsilon_n \varsigma)^\intercal(\sigma+\epsilon_n \varsigma): D^{2}v(t_n,x_n;\sigma)\\
         &-\epsilon_n(\partial_{t}\phi(t_n,x_n)+\frac{1}{2}(\sigma+\epsilon_n \varsigma)^\intercal(\sigma+\epsilon_n \varsigma): D^{2}\phi(t_n,x_n))\\
         &+f(t_n,x_n, \Theta(t_n,x_n,c_n,\epsilon_n))
    \end{aligned}
    \end{equation}
    In the above, for simplifying the notation, we introduce $\Theta(t,x,c,\epsilon):=(v(t,x;\sigma)+\epsilon\phi(t,x)+\epsilon c,\nabla v(,x;\sigma)+\epsilon\nabla\phi(t,x),(\sigma+\epsilon \varsigma)(t,x))$
     Note that by assumption of regularity of $\sigma$, we have 
     \[
     0=-\partial_{t}v(t_n,x_n;\sigma)-\frac{1}{2}\sigma^\intercal\sigma: D^{2}v(t_n,x_n;\sigma)+f(t_n,x_n,\Theta(t_n,x_n,0,0))
     \]
     and, therefore,
      \begin{equation}\label{eqn:apply_regularity}
          \begin{aligned}
         0\ge &\frac{1}{2}\sigma^\intercal\sigma: D^{2}v(t_n,x_n;\sigma)-\frac{1}{2}(\sigma+\epsilon_n \varsigma)^\intercal(\sigma+\epsilon_n \varsigma): D^{2}v(t_n,x_n;\sigma)\\
         &-\epsilon_n(\partial_{t}\phi(t_n,x_n)+\frac{1}{2}(\sigma+\epsilon_n \varsigma)^\intercal(\sigma+\epsilon_n \varsigma): D^{2}\phi(t_n,x_n))\\
         & +f(t_n,x_n \Theta(t_n,x_n,c_n,\epsilon_n)) -f(t_n,x_n, \Theta(t_n,x_n,0,0))
    \end{aligned}
      \end{equation}
     We can use the mean value theorem to obtain the existence of $\bar{\beta}_n$ such that 
     \begin{equation}\label{eqn:mean_value}
         \begin{aligned}
        &f(t_n,x_n, \Theta(t_n,x_n,c_n,\epsilon_n))-f(t_n,x_n, \Theta(t_n,x_n,0,0))\\
         &=\int_0^{1}\dfrac{d}{d\beta}f(t_n,x_n, \Theta(t_n,x_n,c_n,\beta\epsilon_n))d\epsilon\\
         &=\epsilon_n \dfrac{d}{d\beta}f(t_n,x_n, \Theta(t_n,x_n,c_n,\beta\epsilon_n))\Big|_{\beta=\bar\beta_n}
    \end{aligned}
     \end{equation}
     Note that 
     \begin{equation}\label{eqn:expand_d/dbeta}
    \begin{split}
        \dfrac{d}{d\beta}&f(t_n,x_n, \Theta(t_n,x_n,c_n,\beta\epsilon_n))\Big|_{\beta=\bar\beta_n} \\
        =& \epsilon_n\Big(\partial_pf(t_n,x_n, \Theta(t_n,x_n,c_n,\bar\beta\epsilon_n))(\phi(t_n,x_n)+c_n) \\
        &+ \partial_qf(t_n,x_n, \Theta(t_n,x_n,c_n,\bar\beta_n\epsilon_n))\cdot\nabla\phi(t_n,x_n)\\
        &+\partial_{\sigma} f(t_n,x_n, \Theta(t_n,x_n,c_n,\bar\beta_n\epsilon_n)):D^2\phi(t_n,x_n)\Big)
    \end{split}
    \end{equation}
     After using \eqref{eqn:mean_value} and \eqref{eqn:expand_d/dbeta} in \eqref{eqn:apply_regularity} and organizing terms based on the power of $\epsilon_n$, we obtain
     \begin{equation}
         \begin{aligned}
         0\ge \epsilon_n\Big(&-\partial_{t}\phi(t_n,x_n)-\dfrac12\sigma^\intercal\sigma:D^{2}\phi(t_n,x_n)-\sigma^\intercal\varsigma:D^{2}v(t_n,x_n;\sigma)\\
        &+\partial_pf(t_n,x_n,\Theta_n(t_n,x_n,c_n,\bar\beta_n\epsilon_n))(\phi(t_n,x_n)+c_n) \\
        &+ \partial_qf(t_n,x_n,\Theta_n(t_n,x_n,c_n,\bar\beta_n\epsilon_n))\cdot\nabla\phi(t_n,x_n)\\
        &+\partial_{\sigma} f(t_n,x_n,\Theta_n(t_n,x_n,c_n,\bar\beta_n\epsilon_n)):D^2\phi(t_n,x_n)\Big)+o(\epsilon_n)
    \end{aligned}
     \end{equation}
     Dividing both sides by $\epsilon_n$, sending $\epsilon_n\to0$, Assumption \textbf{A3} and using the continuity of $v(\cdot;\sigma)$, $\nabla v(\cdot;\sigma)$, $D^2v(\cdot;\sigma)$, and $\phi$, we obtain
     \begin{equation}
         \begin{aligned}
         0\ge &-\partial_{t}\phi(t,x)-\dfrac12\sigma^\intercal\sigma:D^{2}\phi(t,x)-\sigma^\intercal\varsigma:D^{2}v(t,x;\sigma)\\
        &+\partial_pf(t,x, \Theta(t,x,0,0))\phi(t,x) \\
        &+ \partial_qf(t_n,x_n,\Theta(t,x,0,0))\cdot\nabla\phi(t,x)\\
        &+\partial_{\sigma} f(t,x, \Theta(t,x,0,0)):D^2\phi(t,x)\Big)
    \end{aligned}
     \end{equation}
     with $\Theta(t,x,0,0)=(v(\cdot;\sigma), \nabla v(\cdot;\sigma),\sigma(t,x))$.
     For $(t,x)\in\{T\}\times\mathbb{R}^d$, continuous dependence assumption of $\sigma$ in the direction of $\varsigma$ implies that 
     \begin{equation}
         u^*(T,x)=\limsup_{(s,y,\epsilon)\to(T,x,0)}\frac{v(s,y;\sigma+\epsilon \varsigma)-v(s,y;\sigma)}{\epsilon}\le (1+|x|^p)\lim_{t\to T}\kappa(T-t)=0
     \end{equation}
     which completes the proof.
\end{proof}

\begin{proof}[Proof of Theorem~\ref{thm:max_principle}]
By Lemmas~\ref{lem:u^*_p_growth} and \ref{lem:u^*_sub}, $u*$ and $u_*$ are respectively a sub-solution and a supersolution of \eqref{eqn: directional derivative of u} with with at most $p^\textrm{th}$ growth. Therefore, by 
Assumption \textbf{A3}, regularity of $\sigma$,  $u*\le u_*$. On the other hand, by definition, $u*\ge u_*$. Thus, $u*=u_*$ is the continuous viscosity solution of \eqref{eqn: directional derivative of u}.
\end{proof}
\begin{proof}[Proof of Theorem~\ref{thm:positive_gradient}]
Note that, by Remark~\ref{rem:Girsanov}, Assumption \textbf{A3} ensures that the SDE $dX_s=\mu(s,X_s)ds+\sigma(s,X_s)dB_s$ has a solution. Then, the viscosity solution to equation \eqref{eqn: directional derivative of u} can be given by the Feynman-Kac representation.
  \begin{equation}
      \nabla_\varsigma v(t,x;\sigma)=\mathbb{E}_{t,x}\bigg[\int_t^T e^{-\int_t^s k(r,X_r)dr}\alpha(s,X_s)|\ell|(s,X_s)ds\bigg]
  \end{equation} 
 Finally, $\ell\neq0$ and $\alpha>0$ concludes the proof.
\end{proof}
\begin{proof}
    [Proof of Theorem~\ref{thm:derive_zero}]
    By Theorem~\ref{thm:positive_gradient}, $\ell\equiv 0$, which implies that for all $(t,x)$
    \begin{equation}
    \begin{split}
    \ell(t,x)&=\partial_{\sigma} f\big(t,x,v(t,x;\sigma) ,\nabla v(t,x;\sigma),\sigma(t,x)\big)-\sigma^\intercal D^2v(t,x;\sigma)\\
    &=
    \dfrac{\partial}{\partial \sigma} \Big[
        f\big(t,x,v(t,x;\sigma) ,\nabla v(t,x;\sigma),\sigma(t,x)\big) -\frac12 \sigma^\intercal\sigma : D^2 v(t,x;\sigma)\Big]\!\!=0
    \end{split}
    \end{equation}
    Note that by convexity of $f$ in $\sigma$, this implies that $\sigma$ is the maximizer 
    of 
    \begin{equation}    
    \sigma\mapsto \dfrac12\sigma^\intercal\sigma: D^2v(t,x;\sigma) - f\big(t,x,v(t,x;\sigma),\nabla v(t,x;\sigma),\sigma\big)
    \end{equation}
    Therefore, by \eqref{eqn:conjugate},
      \begin{equation} 
      \begin{split}
          f\big(t,x,v(t,x;\sigma) ,\nabla v(t,x;\sigma),\sigma(t,x)\big)&-\dfrac12(\sigma^\intercal\sigma): D^2v(t,x;\sigma)\\
          &=-h\big(t,x,v(t,x;\sigma),\nabla v(t,x;\sigma),D^2v(t,x;\sigma)\big)
      \end{split}
    \end{equation}
    and 
    \begin{equation} 
      \begin{split}
          0&=\partial_t v(t,x;\sigma) -\dfrac12(\sigma^\intercal\sigma): D^2v(t,x;\sigma)
 + f\big(t,x,v(t,x;\sigma) ,\nabla v(t,x;\sigma),\sigma(t,x)\big)\\
          &=\partial_t v(t,x;\sigma)-h\big(t,x,v(t,x;\sigma),\nabla v(t,x;\sigma),D^2v(t,x;\sigma)\big)
      \end{split}
    \end{equation}
    Finally, Assumption \textbf{A1} implies that $v(t,x;\sigma)$ is the unique viscosity solution of \eqref{eqn:fully_nonlinear}.
\end{proof}
\section{Deep numerical scheme for semilinear and linear equations}
\label{sec:semi_linear_scheme}
In \citet{JHW18}, the connection between a semilinear PDEs and BSDEs  is used to establish a deep scheme for the semilinear PDEs. 
Consider the semilinear PDE 
\begin{equation}\label{app:semi_linear}
 \begin{cases}
 	 	-\partial_t{u}-\frac12(\sigma^\intercal\sigma): D^2 u(t,x)+F(t,x,u(t,x),\nabla u(t,x))=0\\
 	 	u(T,x)=g(x)
 \end{cases}
 \end{equation} 
By the theory of BSDE, \citet{MY99}, the solution of BSDE 
\begin{equation}
\begin{split}
        Y_t&=g(X_T)-\int_t^TF(s,X_s,Y_s,Z_s,\sigma_s)ds-\int_t^T Z_s dX_s\\
    dX_t&=\sigma(t,X_t)dB_t
\end{split}
\end{equation}
is related to the semilinear PDE by $Y_t=u(t,X_t)$ and $Z_t=\nabla u(t,X_t)$. 
In specific case, $F(s,x,p,q)=k(t,x)p+\mu(t,x)\cdot q+\ell(t,x)$ correspond to the linear PDEs.

The main idea is to write the BSDE as a regular forward SDE:
\begin{equation}
\begin{cases}
    dX_t=\sigma(t,X_t)dt\\
    Y_t=Y_0+\int_0^t\big(F(s,X_s,Y_s,Z_s)ds+Z_s\sigma_sdB_s\big)
\end{cases}
\end{equation}
Then, we minimize the loss function
\begin{equation}
    \mathbb{E} \bigg[\Big(\phi(X_0)+\int_0^T\big(F(s,X_s,Y_s,\psi(s,X_s))ds+\psi(s,X_s)\cdot\sigma_sdB_s\big)-g(X_T)\Big)^2\bigg]
\end{equation}
subject to 
\begin{equation}
Y_t=\phi(X_0)+\int_0^t\big(F(s,X_s,Y_s,Z_s)ds+Z_s\sigma_sdB_s\big)
\end{equation}
where the infimum is over all function $\phi$ and $\psi$. If the solution to the PDE, equivalently BSDE, exists, then the risk function vanishes at $\phi^*$ and $\psi^*$ and $Y_0:=\phi^*(X_0)$ and $Z_t:=\psi^*(t,X_t)$ solve the BSDE. Equivalently, $\phi^*(t,x)$ and $\psi^*(t,x)$ equal to $u(t,x)$ and $\nabla u(t,x)$.

\cite{JHW18} approximate $\phi$ and $\psi$ by deep neural networks $Y(x;\theta_0)$ and $Z(t,x;\theta_1)$ and perform the  minimization on $(\theta_0,\theta_1)$: 
\begin{equation}\label{eqn:loss}
    \begin{split}
        \mathop{\inf}\limits_{\theta_0,\theta_1} \mathbb{E}& \bigg[\Big(Y(X_0;\theta_0)+\int_0^T\Big(F(s,X_s,Y_s,Z(s,X_s;\theta_1)\Big)ds\\
        &\hspace{4cm}+Z(s,X_s;\theta_1)\sigma_sdB_s\big)-g(X_T)\Big)^2\bigg]
    \end{split}
\end{equation}
where $Y(s,X_s)$ satisfies
\begin{equation}
Y_{t}=Y(x;\theta_0)+\int_0^t\big(F(s,X_s,Y_s,Z(s,X_s;\theta_1))ds+Z(s,X_s;\theta_1)\sigma_sdB_s\big)
\end{equation}
Practically speaking, after discretization and replacing risk function with an empirical risk function, the minimization yields
\begin{equation}\label{eqn:risk_minimization_semilinear}
        \begin{split}
        (\hat{\theta}_0,\hat{\theta}_1)\in &\underset{\theta \in \mathbb{R}^{k}}{\mathrm{argmin}}\ \mathbb{E}\bigg[\bigg|g(X_N) -Y(X_0;{\theta_0}) -\textstyle\sum_{n=0}^{N-1}\Big(F\big(t_n,X_n,\hat{Y}_{t_n},Z(t_n,X_n;{\theta_1})\big)\Delta t\\
        &\hspace{4cm}-Z(t_n,X_n;\theta_1)\cdot\sigma(t_n,X_n)\Delta B_n\Big) \bigg|^2\bigg]
        \end{split}
    \end{equation}
    subject to 
\begin{equation}\label{app:discrete_Y}
\hat{Y}_{t_n+1}=Y(x;\theta_0)+\textstyle\sum_{i=0}^{n}\big(F(t_i,X_{t_i},\hat{Y}_{t_i},Z(t_i,X_{t_i};\theta_1))\Delta t +Z(t_i,X_{t_i};\theta_1)\sigma_{t_i}\Delta B_{t_{i+1}}\big)
\end{equation}
where $X_0 \sim \mu, \,  
  X_{n+1}:=X_n+\sigma(t_nX_n)\Delta B_n$,  
 $\Delta t = \dfrac{T}{N}$, $t_n=n\Delta n$ and, $\Delta B_n=B_{t_{n+1}}-B_{t_n}$. $Y(x;\hat{\theta}_0)$ and $Z(t,x;\hat{\theta}_1)$ are the approximate solution at time $0$ for the semilinear PDE and the gradient of the solution, respectively. 
\subsection{Modification of the algorithm: global of the value function}
If we need the approximation of the value function $u(t,x)$ for \eqref{app:semi_linear}, one can use \eqref{app:discrete_Y} to generate data for $\hat Y$ and fit a function $\hat u(t,x)$ to the data. However, this adds another computational layer to each iteration of the algorithm. To avoid this, we propose two modifications of \citet{JHW18}. 

First, one can simply take one neural network $U(t,x;\theta)$ to approximate of $u(t,x)$ and replace $Y(x;\theta_0)$ and $Z(t,x;\theta_1)$ by $\nabla U(0,x;\theta)$ and $\nabla U(t,x;\theta)$, respectively. The rest of the method including the loss function does not change. Although this results in similar accuracy, it requires a more complex architecture. 

Second approach is to take two neural network $V(x;\theta_0)$ and $U(t,x;\theta_1)$ as approximations of $u(0,x)$ and  $u(t,x)$, respectively. In the method, we replace $Y(x;\theta_0)$ and $Z(t,x;\theta_1)$ by $V(0,x;\theta_0)$ and $\nabla U(t,x;\theta_1)$, respectively. We also modify the loss function to make sure $U(t,x;\theta_1)$ shifts to the approximation of the value function:
\begin{equation}\label{eqn:risk_minimization_semilinear_mod}
        \begin{split}
        (\hat{\theta}_0,\hat{\theta}_1)\in &\underset{\theta_0,\theta_1}{\mathrm{argmin}}\ \mathbb{E}\bigg[\bigg|g(X_N) -V(X_0;{\theta_0}) -\textstyle\sum_{n=0}^{N-1}\Big(F\big(t_n,X_n,\hat{Y}_{t_n},\nabla U(t_n,X_n;{\theta_1})\big)\Delta t\\
        &-\nabla U(t_n,X_n;\theta_1)\cdot\sigma(t_n,X_n)\Delta B_n\Big) \bigg|^2+|g(X_T)-U(t_0,X_0;{\theta_1})|^2\\
        &\hspace{5cm}+|V(X_0;{\theta_0})-U(t_0,X_0;{\theta_0})|^2\bigg]
        \end{split}
    \end{equation}
    subject to 
\begin{equation}
\hat{Y}_{t_n+1}=V(x;\theta_0)+\textstyle\sum_{i=0}^{n}\big(F(t_i,X_{t_i},\hat{Y}_{t_i},\nabla U(t_i,X_{t_i};\theta_1))\Delta t +\nabla U(t_i,X_{t_i};\theta_1)\sigma_{t_i}\Delta B_{t_{i+1}}\big)
\end{equation}
Out test shows that the second approach suits this study better than the first approach.

\bibliographystyle{plainnat}
\bibliography{refs.bib}
\end{document}